\documentclass{article}
\usepackage{amsmath}
\usepackage{amsthm}
\usepackage{amsfonts}
\usepackage{amssymb}

\usepackage[left=1in,right=1in,top=1in,bottom=1in]{geometry}

\usepackage{hyperref}

\usepackage{fancyhdr}

\numberwithin{equation}{section}

\usepackage{color}
\usepackage{epsfig}

\title{Improved decay rates with small regularity loss for the wave equation about a Schwarzschild black hole}
\author{P. Blue, A. Soffer}

\newtheorem{theorem}{Theorem}[section]

\newtheorem{lemma}[theorem]{Lemma}
\newtheorem{proposition}[theorem]{Proposition}
\newtheorem{definition}[theorem]{Definition}
\newtheorem{remark}[theorem]{Remark}

\newtheorem{faketheoremaux}{Theorem}

\newtheorem{fakecorollaryaux}{Corollary}

\newtheorem{fakelemmaaux}{Lemma}

\newcounter{step}
\newenvironment{stepenv}{}{}
\newcommand{\newstepsequence}{\setcounter{step}{0}}
\newcommand{\newstep}[1]{\addtocounter{step}{1} {\bf Step \thestep: { #1 }} \begin{stepenv}\end{stepenv} }

\newcommand{\alignpartlabelvoodoostart}{\addtocounter{equation}{1}}

\newcommand{\Naturals}{\mathbb N}
\newcommand{\Reals}{\mathbb{R}}

\newcommand{\Fourier}[1]{\mathcal{F}[#1]}
\newcommand{\Lebesgue}{L}

\newcommand{\dalembertian}{\Box}
\newcommand{\defin}{\equiv}
\newcommand{\sgn}{\hbox{sgn}}

\newcommand{\Dt}{\frac{d}{dt}}
\newcommand{\dt}{\partial_t}
\newcommand{\dr}{\partial_{r_*}}
\newcommand{\drr}{\partial_{r_*}^2}
\newcommand{\dtt}{\partial_t^2}
\newcommand{\drrr}{\partial_{r_*}^3}
\newcommand{\dx}{\partial_x}

\newcommand{\phialinsolu}{$\phi\in\solset$ is a solution to the linear wave equation \eqref{eLW}}
\newcommand{\phiasolu}{$\phi\in\solset$ is a solution to the non-linear wave equation \eqref{eNLW}}
\newcommand{\psianicefn}{$\psi\in\solset$}

\newcommand{\starman}{\mathfrak{M}}

\newcommand{\dThree}{d^3\mu}
\newcommand{\dTildeThree}{d^3\tilde{\mu}}

\newcommand{\solset}{\mathbb{S}}

\newcommand{\ConfChgDens}{e_{\mathcal{C}}}
\newcommand{\ConfChg}{E_{\mathcal{C}}}

\newcommand{\linH}{H}

\newcommand{\tl}{\tilde{l}}
\newcommand{\Vtl}{V_l}
\newcommand{\Ptl}{{\bf P}_{l}}
\newcommand{\al}{(\alpha_l)_*}

\newcommand{\chia}{\chi_\alpha}

\newcommand{\Gammahat}{\Gamma}
\newcommand{\CGamma}{C_{\Gammahat}}
\newcommand{\pjdelta}{{ {\boldsymbol \xi}_{j\delta} }}
\newcommand{\phalfminustwodelta}{{ {\boldsymbol \xi}_{\frac12-2\delta} }}
\newcommand{\ErrorTermsMinusTwoDeltaPsi}{O(\| L^{1-2\delta}\chi \psi \|, L^1)}

\newcommand{\Gammacrown}{\bar{\Gamma}}
\newcommand{\chit}{\breve{\chi}}
\newcommand{\chitfn}{\breve{\chi}_{[0,1]}}
\newcommand{\OSterbenz}{O(\chit)}

\newcommand{\pnnan}{\pnna{\pn}}
\newcommand{\pnncn}{\pnnc{\pn}}
\newcommand{\gammam}{\gamma_{L^m}}
\newcommand{\gm}{g_{L^m}}
\newcommand{\gmp}{g_{L^m}'}
\newcommand{\gmpp}{g_{L^m}''}
\newcommand{\gmppp}{g_{L^m}'''}

\newcommand{\supt}{\sup_{0\leq\tau\leq t}}
\newcommand{\ELtwo}{E[L^{2\varepsilon}\phi,L^{2\varepsilon}\dot\phi]}
\newcommand{\ELten}{E[L^{10\varepsilon}\phi,L^{10\varepsilon}\dot\phi]}
\newcommand{\ELeighteen}{E[L^{18\varepsilon}\phi,L^{18\varepsilon}\dot\phi]}
\newcommand{\ELalpha}{E[L^{\alpha}\phi,L^{\alpha}\dot\phi]}
\newcommand{\ConfChgphi}{\ConfChg[\phi,\dot\phi]}

\newcommand{\pniaWITHlambda}[1]{\pniaNOARG(\lambda^{-\delta},#1)}
\newcommand{\pnibWITHlambda}[1]{\pniaNOARG(\lambda^{-\delta},#1)}
\newcommand{\pniWITHlambda}[1]{{\pnifn}_{\lambda^{-\delta}\leq|#1|\leq1}}
\newcommand{\pnplusdelta}{{{ {\boldsymbol \xi}_{n+\delta}}}}

\newcommand{\OSterbenzEpsilon}{O(L^\varepsilon\chit)}

\newcommand{\BSSnorm}[1]{\left\vert\kern-0.9pt\left\vert\kern-0.9pt\left\vert (#1,\dot{#1}) \right\vert\kern-0.9pt\right\vert\kern-0.9pt\right\vert}

\newcommand{\qm}{{\bf x}_m}
\newcommand{\pn}{{ {\boldsymbol \xi}_n }} 

\newcommand{\qn}{{\bf x}_n}

\newcommand{\qmnfn}{X}
\newcommand{\qmn}[1]{\qmnfn_\downarrow(#1)}

\newcommand{\qmnt}[1]{\tilde\qmnfn_\downarrow(#1)}
\newcommand{\qmf}[1]{\qmnfn_\uparrow(#1)}

\newcommand{\ErrorTermsnPsi}{O(\|L^{\frac{1+2n}{2}}\chia \psi\|^2,L^1)}
\newcommand{\bndd}{B}

\newcommand{\pnnfn}{\Phi}
\newcommand{\pnnaNOARG}{\pnnfn_{a,\epsilon}}
\newcommand{\pnna}[1]{\pnnaNOARG(#1)}
\newcommand{\pnnaWITHADERIV}[1]{\pnnaNOARG'(#1)}
\newcommand{\pnnaWITHTWODERIV}[1]{\pnnaNOARG''(#1)}
\newcommand{\pnnbNOARG}{\pnnfn_{b,\epsilon}}
\newcommand{\pnnb}[1]{\pnnbNOARG(#1)}
\newcommand{\pnncNOARG}{\pnnfn_{c,\epsilon}}
\newcommand{\pnnc}[1]{\pnncNOARG(#1)}
\newcommand{\pnn}[1]{{\pnnfn}_{|#1|\leq1}}

\newcommand{\pnifn}{\Psi}
\newcommand{\pniaNOARG}{\pnifn}
\newcommand{\pnia}[1]{\pniaNOARG(L^{-\delta},#1)}

\newcommand{\pnibNOARG}{\pnifn_{1}}

\newcommand{\pnicNOARG}{\pnifn_{2}}
\newcommand{\pnic}[1]{\pnicNOARG(#1)}

\newcommand{\pni}[1]{{\pnifn}_{L^{-\delta}\leq|#1|\leq1}}

\newcommand{\pniWITHl}[1]{{\pnifn}_{l^{-\delta}\leq|#1|\leq1}}

\newcommand{\Gammanm}{\Gamma_{n,m}}
\newcommand{\Gammanhalf}{\Gamma_{n,\frac12}}
\newcommand{\Gammann}{\Gamma_{n,n}}

\newcommand{\ad}{\text{Ad}}
\newcommand{\opnorm}{{}}

\newcommand{\hidenumber}{\nonumber}

\newcommand{\wlap}{\Delta_{\cpctW}}
\newcommand{\cpctW}{\mathfrak{W}}

\begin{document}

\maketitle

\begin{abstract}
  We continue our study of the decoupled wave equation in the exterior of a spherically symmetric, Schwarzschild, black hole. Because null geodesics on the photon sphere orbit the black hole, extra effort must be made to show that the high angular momentum components of a solution decay sufficiently fast, particularly for low regularity initial data. Previous results are rapid decay for regular ($H^3$) initial data \cite{BSterbenz} and slower decay for rough ($H^{1+\epsilon}$) initial data \cite{BlueSoffer3}. Here, we combine those methods to show boundedness of the conformal charge. From this, we conclude that there are bounds for global in time, space-time norms, in particular
\begin{align*}
\int_{\text{I}} |\tilde\phi|^4 d^4\text{vol} <&C , \\
\end{align*}
for $H^{1+\epsilon}$ initial data with additional decay towards infinite and the bifurcation sphere. Here $\tilde\phi$ refers to a solution of the wave equation. $I$ denotes the exterior region of the Schwarzschild solution, which can be expressed in coordinates as $r>2M$, $t\in\Reals, \omega\in S^2$, and $d^4\text{vol}$ is the natural $4$-dimensional volume induced by the Schwarzschild pseudo-metric. 

We also demonstrate that the photon sphere has the same influence on the wave equation as a closed geodesic has on the wave equation on a Riemannian manifold. We demonstrate this similarity by extending our techniques to the wave equation on a class of Riemannian manifolds. Under further assumptions, the space-time estimates are sufficient to prove global bounds for small data, nonlinear wave equations on a class of Riemannian manifolds with closed geodesics. We must use global, space-time integral estimates since $L^\infty$ estimates cannot hold at this level of regularity. 
\end{abstract}


\section{Introduction}

The Schwarzschild solution is a manifold which satisfies the Einstein equation.  Recall that, in general relativity, space-time is described by a four dimensional, Lorentz manifold which satisfies the Einstein equation. In vacuum, this equation reduces to the statement that the Ricci curvature is zero. 

In 1918, Schwarzschild found the first non-trivial solution, which is the unique spherically-symmetric, vacuum solution, other than Minkowski space, flat $\Reals^{3+1}$ with the metric $-dt^2+d\vec{x}^2$. The Lorentz pseudo-metric for the Schwarzschild solution is
\begin{align}
ds^2 = - (1-2M/r) dt^2 + (1 -2M/r)^{-1} dr^2 +r^2 (d\theta^2 + \sin^2(\theta)d\phi^2)
\label{eSchwarzschildMetric}
\end{align}
This is the pseudo-metric outside any spherically symmetric body with mass $M$, and, for large $r$, the space-like geodesic motion on this manifold, as described by the coordinates $(t,r,\theta,\phi)$, approaches that of a body moving under the influence of a central mass, according to Newton's laws, as described by polar coordinates in $\Reals^{3+1}$. Similarly, null geodesics describe the motion of light rays, which physicists treat as massless particles, called photons. 

\begin{figure}
\label{FSchwarzschild}
\begin{center}
\input{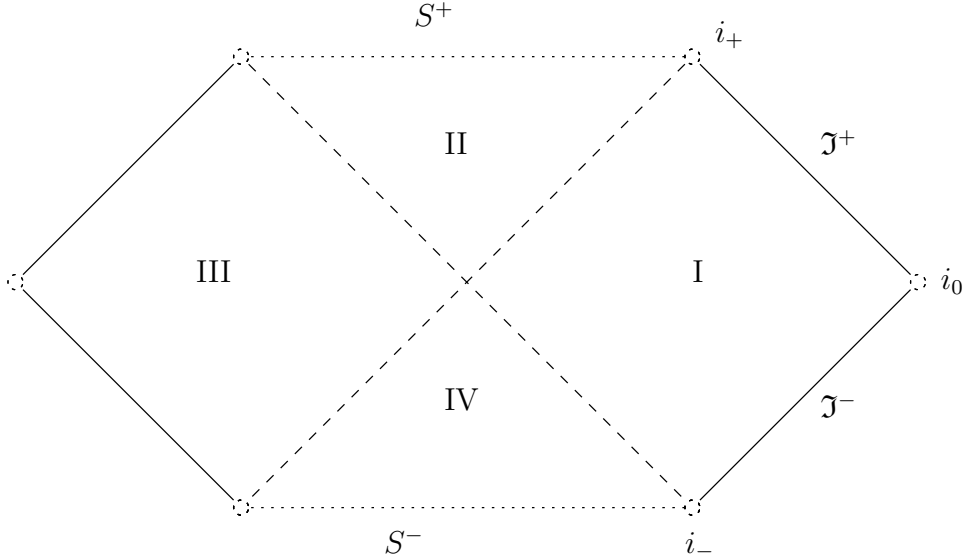}
\end{center}
\caption{The causal structure of the Schwarzschild solution. }
\end{figure}

The pseudo-metric \eqref{eSchwarzschildMetric} is singular at $r=0$ and appears to be singular at $r=2M$, and this was believed to be the case for many years; however, this is an artefact of the coordinate system, and in the completion of the manifold, the surfaces corresponding to $r=2M$ are null geodesic surfaces. The structure of the completion of the Schwarzschild solution is a little surprising, but, by now, well understood \cite{EllisHawking}, and illustrated, up to a conformal transform and ignoring the angular coordinates, in figure \ref{FSchwarzschild}. 

Since as $r\rightarrow\infty$, the Schwarzschild solution approaches Minkowski space, the causal boundary near $r\rightarrow\infty$ is like that in Minkowski space, with a future and past time like infinity, $i^\pm$; a spatial infinity, $i_0$; and future and past null infinities, $\mathfrak{I}^\pm$, which corresponds to a null hyper-surface when the Schwarzschild solution is embedded via a conformal transform ($ds^2\mapsto \Omega^2 ds^2$) in a compact Lorentz space. The intersection of the chronological past of the future null infinity with the chronological future of the past null infinity, $I^-(\mathfrak{I}^+)\cap I^+(\mathfrak{I}^-)$ is commonly referred to as the exterior region of the Schwarzschild black hole and depicted as region $\text{I}$ in figure \ref{FSchwarzschild}. The boundaries of the past of the future null infinity ($I^-(\mathfrak{I}^+)$) and of the future of the past null infinity ($I^+(\mathfrak{I}^-)$) both correspond to surfaces where $r=2M$. The exterior region is described by the Lorentzian pseudo-metric \eqref{eSchwarzschildMetric} on the set 
\begin{align*}
\text{I} = \Reals \times\starman = \{(t,r,\theta,\phi) | t\in\Reals, r>2M, (\theta,\phi)\in S^2 \}
\end{align*}

The singularity at $r=0$ corresponds to both the space like surfaces $S^+$ and $S^-$ in figure \ref{FSchwarzschild}. The regions II and IV around them are both described by the Lorentzian pseudo-metric \eqref{eSchwarzschildMetric} with $r<2M$. Region III is completely isomorphic to region $\text{I}$. The surfaces $r=2M$ are null hyper-surfaces around the region $r<2M$ which are referred to as the black hole, since light rays can enter but not escape from them. More precisely, since events inside the black hole (in region II) can not be witnessed by observers remaining outside (headed to $i^+$ or $\mathfrak{I}^+$), the surfaces $r=2M$ are referred to as event horizons. The coordinate $t$ diverges on the event horizon (diverging to negative infinity in the past, and plus infinity in the future), so that it fails to distinguish points on the event horizon. The points where the four regions meet, corresponding to taking the limit, in region $\text{I}$, as $r\rightarrow2M$ for any fixed value of $t$, is called the bifurcation sphere. 

At $r=3M$, there are null geodesics which orbit the black hole, remaining at $r=3M$. Since photons can orbit here, this surface is known as the photon sphere. Null geodesics initially tangent to a surface of constant $r$ will escape to infinity if $r>3M$ and fall into the event horizon, and then into the singularity, if $r<3M$. From this perspective, $r=3M$ is a natural boundary. However, this surface does not play a special role in terms of causal geometry. 

A major open problem in relativity is the stability of black holes under small perturbations. This is important for understanding the behavior, interaction, and dynamics of black holes. Exact solutions exist for black holes with electromagnetic charge and angular momentum, as well as mass. Charged solution are given by the Reissner- Nordstr{\o}m solutions. Black holes with angular momentum rotate and are described by Kerr solutions, or Kerr-Newman solutions if they are also charged. If the charge and angular momentum are not excessive, an event horizon and causal infinities still exist, and the geometry of region $\text{I}$ is not greatly altered; although, the nature of the singularities changes dramatically and there are infinitely many additional exterior regions, like region III in the Schwarzschild solution. By the ``no hair theorem'', stationary black holes are entirely characterized by their mass, charge, and angular momentum, and it is generally believed that black holes approach one of the known exact solutions. From Birkhoff's theorem, it is known that there can be no spherically-symmetric, vacuum (or vacuum-Maxwell) perturbations \cite{EllisHawking}. Since the Einstein equation is nonlinear, it is expected to couple higher spherical harmonic perturbations so that higher spherical harmonic perturbations to the Schwarzschild solution will introduce rotation, and that, although the Kerr solutions are expected to be stable as a class, the Schwarzschild solution alone will not be. However, since the Kerr solution has no global time-like Killing vector (even just in the exterior), its perturbations are extremely hard to analyze. Perturbation equations for Schwarzschild and Kerr-Newman have been derived by Regge and Wheeler \cite{ReggeWheeler}, Zirilli \cite{Zirilli}, and Teukolsky \cite{Teukolsky}, and all of these suggest black holes are stable. Whiting has shown there are also no unstable modes \cite{Whiting}. 

The only rigorous proof of nonlinear stability is that of Christodoulou and Klainerman, which proves the stability of Minkowski space \cite{ChristodoulouKlainerman} (see also \cite{LindbladRodnianski} for a simplified proof). All the perturbation equations derived and the analysis of Minkowski stability have involved second order wave equations, so that we hope that information on the decay of the linear wave equation will provide information for the stability of Kerr solutions. 

The rigorous analysis of the wave equation on the Schwarzschild solution began with \cite{DimockKay} which proved the $L^\infty$ boundedness of solutions with initial data in $H^2$. Scattering results have been proven for wave equations on the Schwarzschild \cite{BachelotNicolas} and Kerr \cite{Hafner} solutions. Decay for various types of fields, including the wave equation, around a Kerr black hole have been proven using spectral techniques \cite{FinsterKamranSmollerYau}. A similar approach to ours has been used to prove decay for the wave equation on the Schwarzschild solution, including estimates along the event horizon, known as Price's law\cite{DafermosRodnianski}. Estimates along the event horizon are important not only for perturbation in region $\text{I}$, but also for understanding the nature of the singularity \cite{Dafermos}

We consider the decoupled wave equation on the exterior of the Schwarzschild black hole, 
\begin{align}
\label{eLWTrueInIntro}
\dalembertian \tilde\phi =& 0 ,& \tilde\phi(0)=&\tilde\phi_0,&\dt\tilde\phi(0)=&\tilde\phi_1 ,
\end{align}
where the d'Alembertian $\dalembertian$ is given by the fixed background pseudo-metric \eqref{eSchwarzschildMetric}. By introducing the Regge-Wheeler tortoise coordinate, $r_*$, and transforming to $\phi=r\tilde\phi$ (and similarly for the initial data), we transform the wave equation to 
\begin{align}
\label{eLWInIntro}
-\dtt\phi=&-\drr\phi +V\phi +V_L(-\Delta_{S^2})\phi , &
\phi(0)=&\phi_0, &\dt\phi(0)=&\phi_1 .
\end{align}
This has the form of a family of one dimensional wave equations on each spherical harmonic, and we have previously treated it as such \cite{BlueSoffer3,BSterbenz}. Our method is to introduce the energy and conformal charge, 
\begin{align}
\label{eEnergyConfChgInIntro}
E=& \int_{\{t\}\times\starman} |\dt\phi|^2+ V|\phi|^2 + V_L |\nabla_{S^2}\phi|^2 dr_*d^2\omega,\\
\ConfChg=& \int_{\{t\}\times\starman} (t^2+r_*^2)(|\dt\phi|^2+ V|\phi|^2 + V_L |\nabla_{S^2}\phi|^2) + 2tr_*(dt\phi)(\dr\phi) dr_* d^2\omega + E .\nonumber
\end{align}
We analyze this using commutator methods and a radial derivative operator. Since the energy comes from a Killing vector, it has long been known to be conserved. The conformal charge has been used to analyze wave equations in $\Reals^{n+1}$, for example \cite{GinibreVelo}. Since the conformal charge, on the Schwarzschild solution, does not come from a Killing vector, it is not conserved. 

Previously, we have shown that the growth of the conformal charge can be estimated by the accumulation of the angular component of the energy which does not disperse from near the photon sphere, 
\begin{align}
\label{eConfChgEstimateInIntro}
\ConfChg(t_1)-\ConfChg(t_0) \leq& \int_{t_0}^{t_1}\int_\starman \chi (|\phi|^2 + |\nabla_{S^2}\phi|^2) dr_* d^2\omega dt .
\end{align}
Here $\chi$ is a compactly supported function coming from the positive part of the trapping terms for the potentials, 
\begin{align*}
2V+r_*V', 2V_L+r_*V_L' .
\end{align*}
Our analysis begins with a Morawetz estimate\footnote{
Morawetz originally used a somewhat different estimate, which was modified to prove local decay estimates of the type we use in $\Reals^{n}$ \cite{Lavine} and later for radial data on the Schwarzschild solution \cite{LabaSoffer}. Morawetz also introduced what we refer to as the conformal estimate in $\Reals^n$, which some authors refer to as the Morawetz estimate. } 
\begin{align}
\label{eLocDecInIntro}
\int_0^\infty \int_\starman \chi (|\phi|^2 + |\dr\phi|^2)  dr_* d^2\omega dt \leq C E .
\end{align}
This estimate is proven using the technique first introduced in \cite{Morawetz} but focuses on the role of different terms. Since angular derivatives can be applied to the linear wave equation \eqref{eLWInIntro}, we can apply the Morawetz estimate \eqref{eLocDecInIntro} to \eqref{eConfChgEstimateInIntro} to get that the conformal charge is bounded in terms of $\phi$ and the operator square root of $(1-\Delta_{S^2})$, denoted $L$, by
\begin{align}
\label{ePrelimEstimateInIntro}
\ConfChg(t)\leq& Ct E[L\phi] + \ConfChg(0). 
\end{align}
Since $E$ is an $L^2$ norm of derivatives of $\phi$, and $L$ acts like an additional derivative $\nabla_{S^2}$, we refer to \eqref{ePrelimEstimateInIntro} as an $H^2$ estimate. 

We are interested in improving \eqref{ePrelimEstimateInIntro} since the growth in $t$ of the conformal charge does not give good decay estimates for $\phi$ or $\tilde\phi$ and the loss of angular regularity can make it hard to close nonlinear estimates. In \cite{BSterbenz}, it was possible to improve \eqref{ePrelimEstimateInIntro} to uniform boundedness in time  for the conformal charge by proving a Morawetz estimate localized inside the light cone and using $t^{-2}\ConfChg$ instead of $E$ as the bound in the Morawetz estimate \eqref{eLocDecInIntro}. To extend to non-linear problems, $L^\infty$ estimates were used, which we required an additional angular derivative, leaving a $H^3$ estimate. In \cite{BlueSoffer3}, we were able to reduce the loss of angular regularity in \eqref{eLocDecInIntro} from one derivative to $\varepsilon$ derivatives, which is defined by the spectral theorem. We did this by dilating the radial derivative operator by fractional powers of the angular derivative operator and by relating powers of the angular derivative to the radial derivative, techniques we refer to as angular modulation and phase space analysis. Here, we combine these techniques. Since we are working only at an $H^{1+\varepsilon}$ level of regularity, it is not possible to prove an $L^\infty$ estimate, and, instead we prove global-in-time, space-time integral estimates, the most transparent of which is

\begin{proposition} \label{PSchwarzschildSTNorm} If $\tilde\phi$ is a solution of \eqref{eLWTrueInIntro} and $\varepsilon>0$, then 
\begin{align*}
\int_{\text{Region $\text{I}$}} |\tilde\phi|^4 d^4\text{vol} 
<C (\ConfChg[\phi_0,\phi_1] + E[\phi_0,\phi_1] + E[L^\varepsilon\phi_0,L^\varepsilon\phi_1])
\end{align*}
provided the norms on the right are finite, the initial data vanishes sufficiently fast at the bifurcation sphere and spatial infinity, and where $\phi=r\tilde\phi$ and $d^4\text{vol}$ is the four dimensional space-time volume associated to the pseudo-metric in \eqref{eSchwarzschildMetric}. 
\end{proposition}

The natural measure $d^4\text{vol}$ can be written in Schwarzschild coordinates as $\sqrt{g}dx^1\ldots dx^4 = r^2\sin(\theta) dr d\theta d\phi dt$. 

The estimate in \eqref{eConfChgEstimateInIntro} reveals the importance of the photon sphere, $r=3M$, in this problem. The null geodesics which orbit at $r=3M$ appear to have some similarity to closed geodesics on a Riemannian manifold, since their spatial coordinates are periodic in time\footnote{We do not claim that the behavior of null geodesics tangent to the photon sphere is analogous to closed null geodesics on a Lorentz manifold, which violate causality and which are not present in the Schwarzschild manifold. }. However, because of the curvature of both space and time, it is not clear how closely this analogy holds. Physicists are well aware of various strange effects caused by the photon sphere, for instance, an infinite number of visible images are created by gravitational lensing around the photon sphere \cite{GravLensing}. 

The theory of dispersive waves on manifolds is still developing, but dispersion is known to be weaker in the presence of closed geodesics. Although the exponents differ, in $\Reals^{n+1}$, the wave and Schr\"odinger equations have similar $L^1\rightarrow L^\infty$ pointwise in time decay rates and similar Strichartz estimates \cite{Cazenave,Sogge}. In the absence of closed geodesics, global Strichartz estimates have been proven for both the wave \cite{SmithSogge} and Schr{\"o}dinger \cite{HassellTaoWunsch}. For the wave equation, because of finite speed of propagation, the local in time theory is the same as in $\Reals^{n+1}$; however, for the Schr{\"o}dinger equation on a compact manifold with closed geodesics, on a finite time interval, Strichartz estimates hold, but with an $\varepsilon$ loss of regularity \cite{BurqGerardTzvetkov}. Similarly, while an estimate like \eqref{eLocDecInIntro} holds for the Schr{\"o}dinger equation in $\Reals^{n+1}$ (where it is known as a local smoothing estimate, since control on the left is by an $H^{1/2}$ norm for the Schr{\"o}dinger equation), there must be an additional $\varepsilon$ derivative loss in the presence of closed geodesics \cite{Doi}. For the wave equation in the presence of closed geodesics, it is not possible for an estimate like \eqref{eLocDecInIntro} to hold if it controls derivatives in all directions, since, through a geometric optics approximation, wave equation solutions in the energy space can be shown to stay arbitrarily close to the path of geodesics \cite{Ralston}. Since we are working with commutator methods, which typically can also be applied to the Schr{\"o}dinger equation, we do not expect to be able to avoid an $\varepsilon$ regularity loss in the presence of closed geodesics. 

The analogy we wish to consider is to a wave equation on a three-dimensional, warped product manifold, $\starman=\Reals\times\cpctW$, with metric $ds^2=dr_*^2 + r(r_*)^2 d\omega^2$ and $\cpctW$ compact. We make several assumptions on $r(r_*)$, explained in section \ref{sSetup}, but essentially that $r$ has a unique minimum, which corresponds to a single closed geodesic surface, and that $r$ grows at least linearly but at most polynomially. This Riemannian metric defines a Laplace- Beltrami operator, which we can use to define a wave equation, 
\begin{align}
\label{eWarpedProductWaveInIntro}
-\dtt\tilde\phi =& -\Delta_{\starman}\tilde\phi ,& \tilde\phi(0)=&\tilde\phi_0,&\dt\tilde\phi(0)=&\tilde\phi_1.
\end{align}
Introducing $\phi=r\tilde\phi$ again, we get an equation of the form \eqref{eLWInIntro}. The analysis we develop for the wave equation on the Schwarzschild solution applies equally well to this problem both strengthening the analogy between the photon sphere and a closed geodesic on a Riemannian manifold and proving new estimates for solutions to wave equations in the presence of closed geodesics\footnote{As noted in \cite{BSterbenz}, the higher regularity analysis there also applies to warped product manifolds with a single closed geodesics under certain conditions; although, the necessary conditions were not given. }. We are able to prove the analog of proposition \ref{PSchwarzschildSTNorm}. In fact, for the defocusing, semi-linear wave equation with exponent $8/3 < p < 3$, 
\begin{align}
\label{eNLWarpedProductWaveInIntro}
-\dtt\tilde\phi =& -\Delta_{\starman}\tilde\phi + |\tilde\phi|^{p-1}\tilde\phi ,& \tilde\phi(0)=&\tilde\phi_0,&\dt\tilde\phi(0)=&\tilde\phi_1.
\end{align}
our space-time integral estimates are sufficiently strong to close a boot strap argument and extend our result to the small-data, semi-linear case. 

\begin{proposition} \label{PWarpedProductSTNorm}
If $\tilde\phi$ is a solution to the linear wave equation \eqref{eWarpedProductWaveInIntro} on a warped product manifold, $\varepsilon>0$, and conditions \ref{CPtRadial} to \ref{CNLPotLTheSame} given in section \ref{sSetup} hold, then 
\begin{align*}
\int_{-\infty}^{\infty} \int_\starman |\tilde\phi|^4 d^3\text{vol} dt 
< (\ConfChg[\phi_0,\phi_1] + E[L^\varepsilon\phi_0,L^\varepsilon\phi_1])
\end{align*}
provided that the norms on the left are sufficiently small. The operator $L$ denotes the operator square of $(1-\wlap)$. For $8/3<p<3$ and sufficiently small initial data, in the norms appearing on the right hand side, the same estimate is true if $\tilde\phi$ is a solution to the semi-linear wave equation \eqref{eNLWarpedProductWaveInIntro}. 
\end{proposition}

Because this result replies upon commutator methods, it can also be extended to the Schr\"odinger equation on a warped product manifold, using the methods in this paper and in \cite{BSterbenz} and \cite{LabaSoffer}. We will deal with this elsewhere. 

We begin in section \ref{sSetup} by listing the conditions on the terms in the wave equation \eqref{eLWInIntro} which allow us to analyze it and explaining how to transform the wave equation on the Schwarzschild solution or a warp product manifold to fit this description. In section \ref{sMethodsPD}, we define the energy and conformal charge. In section \ref{sThreeEstimates}, we prove a Sobolev estimate, an estimate on weighted $\Lebesgue^q$ norms in terms of the conformal charge, and a Hardy type estimate. Following this, in section \ref{sHeisenberg}, we introduce an analog of the Heisenberg relation for the wave equation, which we use in section \ref{sLocalDecay} to prove the local decay estimate \eqref{eLocDecInIntro}. In section \ref{sAngularRefinements}, we define and outline the use of the phase space observable, $\Gammahat$, which allows us to improve the local decay estimate by a factor of $L^{1-\varepsilon}$. In section \ref{sTemporalRefinements}, we then localize this inside the light-cone to get additional decay in time, which allows us to close the estimates and prove Propositions \ref{PSchwarzschildSTNorm} and \ref{PWarpedProductSTNorm}. In the appendix, we present a brief version of the arguments in \cite{BlueSoffer3} which allows us to estimate the partial phase space observables, $\Gammanm$, which make up $\Gammahat$.

\section{Setup}
\label{sSetup}

In this section, we state the form of a wave equation to which the rest of our analysis will apply, and we show that the wave equation on both a warp product manifold satisfying certain conditions and on the Schwarzschild solution can be transformed to this form. In addition to the differential equation, we will also introduce function spaces. 

We will consider a wave equation of the form
\begin{align}
\label{eLW}
-\dtt \phi=& -\drr\phi +V\phi +V_L(-\wlap)\phi , &
\phi(1)=&\phi_0, &
\dot\phi(1)=& \phi_1 .
\end{align}
for $\phi(t,r_*\omega)$ with $(t,r_*,\omega)\in\Reals\times\starman=\Reals\times\Reals\times\cpctW$ and $\cpctW$ a two dimensional, compact, smooth manifold. On $\cpctW$, the Laplace- Beltrami operator, $(-\wlap)$, has a spectrum, $\sigma(-\wlap)$, which is discrete \cite{Berard}. Following the standard terminology on the sphere, we refer to the eigenspaces as harmonics. We index the eigenvalues in increasing order by $l$, and denote the projection onto the corresponding eigenspace by $\Ptl$ and the eigenvalue by $\tl^2$. Thus, on the sphere, we have $\tl^2=l(l+1)$. We will also consider the nonlinear wave equation, 
\begin{align}
\label{eNLW}
-\dtt \phi=& -\drr\phi +V\phi +V_L(-\wlap)\phi + f(r_*) |\phi|^{p-1}\phi, &
\phi(1)=&\phi_0, &
\dot\phi(1)=& \phi_1 .
\end{align}
We denote the nonlinearity $F'(|\phi|^2)\phi$, where 
\begin{align}
F(|\phi|^2)=\frac{|\phi|^{p+1}}{p+1}\hidenumber 
\end{align}
We treat the linear wave equation \eqref{eLW} as a special case of the nonlinear wave equation \eqref{eNLW} with $f=0$. 

Since these equations are extremely general, we only consider the case where the following assumptions hold: 
\begin{enumerate}
\item \label{CPtRadial} The potentials, $V$ and $V_L$, are smooth functions and depend only on the radial coordinate, $r_*$. 
\item \label{CPtPos} The potentials are non-negative: $V\geq0$, $V_L\geq0$. 
\item \label{CPtRepulsive} For each $\tl^2$ in the spectrum of $(-\wlap)$, for the effective potential $\Vtl=V+\tl^2V_L$, 
\begin{enumerate}
\item \label{CPtOnePeak} $\Vtl$ has a unique critical point, which is a maximum and occurs at $r_*=\al$, 
\item \label{CPtVanishLinear} $\Vtl$ vanishes linearly at $r_*=\al$, 
\item \label{CPtOriginalAlso} the potentials $V$ and $V_L$ also satisfy conditions \ref{CPtOnePeak} and \ref{CPtVanishLinear} with $(\alpha_0)_*$ and $(\alpha_\infty)_*$ respectively. 
\item \label{CPtalConverge} the $\al$ converge to $(\alpha_\infty)_*$. 
\end{enumerate}
\item \label{CPtCptTrap} The trapping terms, $2V+r_*V'$ and $2V_L+r_*V_L'$, are positive only in a compact region. 
\item \label{CNLrNice} The weight $V_L'/V_L$ is bounded above by a constant, and the angular potential decays at least like $r_*^{-2}$: $V_L'/V_L < C$, $V_L< C(1+r_*^2)^{-1}$.  
\end{enumerate}
For non-linear problems, we require the following additional assumptions:  
\begin{enumerate}
\setcounter{enumi}{5}
\item \label{CNLrBetter} The weight $V_L'/V_L$ decays like $r_*^{-1}$, and the angular potential decays at least like $r_*^{-2}$: $V_L'/V_L < C(1+|r_*|)^{-1}$, $V_L< C(1+r_*^2)^{-1}$. 
\item \label{CNLCoeffPos} The non-linear coefficient is a smooth function which depends only on the radial coordinate and is nonnegative: $f\geq0$. 
\item \label{CNLCptTrap} The non-linear trapping term, $2f+r_*f'$, is positive only in a compact region on which $f$ is strictly positive.  
\item \label{CNLPotLTheSame} The quantity $f V_L^{(1-p)/2}$ is bounded. 
\end{enumerate}
Note that condition \ref{CNLrNice} follows immediately from condition \ref{CNLrBetter}. Note also that the existence of a unique critical point which is a maximum, from condition \ref{CPtOnePeak}, guarantees that the potentials are bounded. Furthermore, from condition \ref{CPtCptTrap} on the trapping term, the quantity $r_*^2V$ and $r_*^2V_L$ are bounded, so the potentials must vanish rapidly as $r_*\rightarrow\pm\infty$. For the local decay, angular modulation, and phase space analysis techniques used in sections \ref{sLocalDecay} and \ref{sAngularRefinements}, conditions \ref{CPtRadial} - \ref{CPtalConverge} are sufficient. 

We will work with the measure $\dThree=dr_*d\omega$. Unless otherwise specified, $L^q$ norms are with respect to this measure, $\|\psi\|$ refers to the $L^2$ norm, and inner products are taken in $L^2$. There is a Fourier transform on $\starman=\Reals\times\cpctW$ defined by decomposing in the eigenbasis for $(-\wlap)$ and then applying the one dimensional Fourier transform on each eigenspace. The Schwartz class, $\solset$, is defined to be the set of pairs $(\psi_0,\psi_1)$ for which
\begin{align}
\tag{S}
\label{SchwartzClass}
\forall i,j,k\in\Naturals, a\in\{0,1\}: \exists C_{a,i,j,k}: \forall (r_*\omega)\in\starman: |r_*^i \dr^j \nabla_{\cpctW}^k \psi_{a}(r_*,\omega) | < C_{a,i,j,k} .
\end{align}
Abusing notation, we will use $\psi(t)\in\solset$ to mean $(\psi(t),\dot\psi(t))\in\solset$. In general, we will use $\phi$ to denote a solution to \eqref{eNLW}, $\psi$ or $\psi(t)$ to denote a function in $\solset$ in the abused notation, and $(v,w)$ to denote a pair in $\solset$. 

In our analysis, the function $\phi$ will be related to a function, $\tilde\phi$, which is the solution to a geometric equation. These are related by $\phi=r(r_*) \tilde\phi$ for some weight function $r(r_*)$. 

We now consider the wave equation on a warp product manifold. We take the manifold to be $(r_*,\omega)\in\starman=\Reals\times\cpctW$, where $\cpctW$ is a two dimensional, compact, smooth, Riemannian manifold with the metric $d\omega^2$. We take the metric on $\starman$ to be 
\begin{align}
ds^2 =& dr_*^2 + r(r_*)^2 d\omega^2 .
\hidenumber \end{align}
This defines a Laplace- Beltrami operator on $\starman$, and we consider the nonlinear wave equation for $\tilde\phi$, 
\begin{align}
-\dtt\tilde\phi=& -\Delta_{\starman}\tilde\phi + |\tilde\phi|^{p-1}\tilde\phi .
\hidenumber \end{align}
Under the transformation $\phi=r\tilde\phi$, this becomes
\begin{align}
-\dtt\phi=& -\drr\phi + V\phi + V_L(-\wlap)\phi + fF'(|\phi|^2)\phi
\hidenumber \end{align}
with
\begin{align}
V=& r'' r^{-1}, \hidenumber\\
V_L=& r^{-2}, \hidenumber\\
f=& r^{1-p}, \hidenumber\\
F(|\phi|^2)=&\frac{1}{p+1} |\phi|^{p+1} .
\hidenumber \end{align}
This is of the form \eqref{eNLW}. If $r$ and $r''$ are positive then conditions \ref{CPtRadial}, \ref{CPtPos}, \ref{CNLCoeffPos}, and \ref{CNLPotLTheSame} are clearly satisfied. Conditions \ref{CPtRepulsive}, \ref{CPtCptTrap}, \ref{CNLrBetter}, and \ref{CNLCptTrap} are, therefore, the only restrictions. From the example $r=(1+r_*^2)$, it is clear these conditions permit some manifolds with closed geodesics. If we consider the restriction on $V_L$ to be the central ones, then condition \ref{CPtRepulsive} says there there is a unique, closed geodesic surface which is unstable. If $r$ behaves polynomially, with $r=r_*^{p_1} + O(r_*^{p_1-1})$ and $r'=r_*^{p_1-1} + O(r_*^{p_1-2})$, then conditions \ref{CPtCptTrap} and \ref{CNLrBetter} follow if $r$ grows faster than linearly, $p_1>1$ (and will also hold if $p_1=1$ provided the lower order terms give condition \ref{CPtCptTrap}). This is also sufficient for \ref{CNLCptTrap} to hold. We work with the measure $dr_*d^2\omega$ for $\phi$ and with $r^2dr_*d^2\omega$ for $\tilde\phi$. This gives an isometry of $L^2$ spaces. 

We now turn to the wave equation on the Schwarzschild manifold. Since there is curvature in time as well as space, the wave equation does not arise from a warp product of the form previously considered. 

The exterior region of the Schwarzschild solution is given by $(t,r,\omega)\in\Reals\times(2M,\infty)\times S^2$ with the metric
\begin{align}
ds^2=& -(1-2M/r)dt^2 +(1-2M/r)^{-1}dr^2 + r^2d\omega^2 .
\hidenumber \end{align}
We introduce the Regge-Wheller tortoise coordinate, $r_*$, defined by 
\begin{align}
\frac{\text{d}r}{\text{d}r_*} =& 1-2M/r ,& r(0)=& 3M . \hidenumber 
\end{align}
It is sometimes useful to work with the explicit expansion, 
\begin{align}
r_* = r +2M \log\left(\frac{r-2M}{2M}\right) -3M +2M\log 2 .
\hidenumber \end{align}
This transforms the exterior region to $(t,r_*,\omega)\in\Reals\times\starman=\Reals\times\Reals\times S^2$. 

The Lorentzian metric defines a d'Alembertian operator and thus the wave equation \eqref{eLWTrueInIntro} for $\tilde\phi$. Using the $r_*$ coordinate and the transformation $\phi=r\tilde\phi$, the wave equation takes the form \eqref{eLW} with
\begin{align}
V=& \frac{2M}{r^3}(1-2M/r) ,\hidenumber\\
V_L=& \frac{1}{r^2} (1-2M/r) , \hidenumber
\hidenumber \end{align}
Once again, it is clear that conditions \ref{CPtRadial}, \ref{CPtPos}, and \ref{CNLrNice}
are satisfied, since we are working in the region $r>2M$. We show conditions \ref{CPtRepulsive} and \ref{CPtCptTrap} hold below. We use the measure $dr_*d^2\omega$ for $\phi$ and use the natural measure $r^2 dr d^2\omega dt=\sqrt{-g}dr d^2\omega dt=(1-2M/r)r^2dr_*d^2\omega dt$ when integrating $\tilde\phi$ in space-time. 

If we were to introduce a semi-linearity into the wave equation, this would translate into a term $fF'(|\phi|^2)\phi=(1-2M/r) r^{1-p} |\phi|^{p-1} \phi$ in equation \eqref{eNLW}. This clearly satisfies conditions \ref{CNLCoeffPos} and \ref{CNLPotLTheSame}. For $p\geq3$, it is not hard to show \ref{CNLCptTrap}. Unfortunately, condition \ref{CNLrBetter} does not hold towards the event horizon, $r\rightarrow2M$, which prevents us from getting a good bound on $\|r^{-1}\phi\|$ in the Sobolev estimate and from controlling the growth of higher Sobolev norms. 


We first show that on each spherical harmonic, the effective potential has a unique maximum. In \cite{BlueSoffer3}, we showed the same holds true on the more general Reissner-Nordstr\o m manifold, but, for simplicity, we restrict to the Schwarzschild manifold here. 

\begin{lemma}
\label{lVtlUniquePeak}
For $l\geq0$ and $M>0$, the effective potential $\Vtl$ has a unique critical point, which is a maximum. The maxima, $\al$, converge to $0$. 
\end{lemma}
\begin{proof}
From the effective potential, we can compute its derivative, which has two roots, $\alpha_{l,\pm}$.
\begin{align}
\Vtl=& V + \tl^2 V_L, \hidenumber\\
\Vtl'=&-2r^{-5}( \tl^2r^2 -(\tl^2-1)3Mr-8M^2)(1-2M/r), \label{eVtlDeriv}\\
\alpha_{l,\pm}=& \frac{ (\tl^2-1)3M \pm\sqrt{(\tl^2-1)^29M^2+32\tl^2M^2}}{2\tl^2} .
\hidenumber \end{align}
For $\tl=0$, the roots are $0$ and $8M/3$. For $\tl>0$, $\alpha_{l,-}$ is negative, and $\alpha_l=\alpha_{l,+}$ is larger than $2M$. At $r=\alpha_l$, $\Vtl'$ goes from negative to positive. From \eqref{eVtlDeriv}, it is clear that $\alpha_l\rightarrow 3M$. The corresponding values of $r_*$ are denoted $\al$ and converge to $(\alpha_\infty)_*=0$. 
\end{proof}

By direct computation, we now show that the trapping terms $2V+r_*V'$ and $2V_L+r_*V_L'$ are positive only in a compact region. 
\begin{align}
2V+r_*V'
=&2(1-\frac{2M}{r})\frac{2M}{r^3}+r_*(1-\frac{2M}{r})\frac{2M}{r^5}(-3r+8M)\hidenumber\\
=& \frac{2M}{r^5}(1-\frac{2M}{r})(2r^2-r_*(3r-8M))\hidenumber\\
2V_L+r_*V_L'
=&\frac{2}{r^2}(1-\frac{2M}{r})+r_* \frac{1}{r^4}(-2r+6M)\hidenumber\\
=&\frac{2}{r^4}(1-\frac{2M}{r})(r^2-r_*(r-3M))
\hidenumber \end{align}
Since for large $r$, $r_*\sim r + \log r$, the terms $2r^2-r_*(3r-8M)\sim -r^2$ and $r^2-r_*(r+3M) \sim -r\log r$ both diverge to negative infinite. For $r_*$ going to negative infinite, $3r-8M<0$ and $r-3M<0$, so that $-r_*(3r-8M)$ and $-r_*(r-3M)$ both diverge to negative infinite like a constant times $r_*$. This proves that the trapping terms $2V+r_*V'$ and $2V_L+r_*V_L'$ are positive only in a compact region. This shows condition \ref{CPtCptTrap} holds. 

From condition \ref{SchwartzClass}, Schwartz class functions and all their derivatives must vanish at the bifurcation sphere. The normalized radial derivative is $(1-2M/r)^{1/2}\dr \psi$. Since $r-2M$ and $r_*$ are related exponentially as $r\rightarrow 2M$, the polynomial decay rate for $r_*$ derivatives does not imply decay, or even boundedness, of $r$ derivatives or normalized radial derivatives, which more accurately describe the behavior of $\psi$ as $r\rightarrow2M$. Thus, functions like $(r-2M)^{\epsilon}$ satisfy condition \ref{SchwartzClass} as $r_*\rightarrow-\infty$ (some other factor which decays as $r\rightarrow\infty$ can then be applied to enforce the condition in the other direction).

\section{Methods for pointwise decay}
\label{sMethodsPD} 

In this section, we summarize the use of various densities, in particular the energy density, which we have previously used to study the wave equation on the Schwarzschild space \cite{BlueSoffer3}. This is an application of very commonly used techniques. 

The energy, radial momentum, and ($T\cpctW$ valued) angular momentum densities are defined by 
\begin{align}
e[v,w]=&\frac12( w^2+v'^2+Vv^2+V_L \nabla_{\cpctW}v\cdot \nabla_{\cpctW}v + f(r_*) F(v^2)) ,\hidenumber\\
p_*[v,w]=& wv' ,\hidenumber\\
\vec{p}_\omega[v,w]=& w\nabla_{\cpctW}v .
\hidenumber \end{align}
If $\phi$ is a solution to the wave equation, these satisfy
\begin{align}
\alignpartlabelvoodoostart
\label{e11.2a}
\tag{\theequation a}
0=& 
\dt e[\phi,\dot{\phi}] 
-\dr p_*[\phi,\dot{\phi}] 
-\nabla_{\cpctW}\cdot\vec{p}_\omega[\phi,\dot{\phi}] , \hidenumber\\
\label{e11.2b} 
\tag{\theequation b}
0=& 
\dt p_*[\phi,\dot{\phi}] 
-\frac12 \dr[ \dot{\phi}\dot{\phi} -\phi'\phi' -\phi V\phi -\nabla_{\cpctW}\phi\cdot V_L \nabla_{\cpctW}\phi -f(r_*)F(\phi^2) ]\hidenumber\\ 
&-\nabla_{\cpctW}\cdot(\phi'V_L\nabla_{\cpctW}\phi)\notag\hidenumber\\
&-\frac12 \phi V'\phi - \frac12 \nabla_{\cpctW}\phi\cdot V_L'\nabla_{\cpctW}\phi -\frac12 f'(r_*) F(|\phi|^2)  .\notag
\hidenumber \end{align}
From integrating equation \eqref{e11.2a}, we have that the total energy, 
\begin{align}
E[\phi,\dot{\phi}] = \int_\starman e[\phi,\dot{\phi}] \dThree 
= \int_{\starman} \frac12( \dot{\phi}^2+\phi'^2+V\phi^2+V_L \nabla_{\cpctW}\phi\cdot \nabla_{\cpctW}\phi + f(r_*)F(\phi^2)) \dThree
\hidenumber \end{align}
is conserved. Since the energy density dominates the $L^2$ norm of the time and radial derivatives of $\phi$, we can use many of the basic arguments which were developed in $\Reals^{n+1}$; however, since the weight, $V_L$, may have more complicated dependency, we will only be able to use the angular derivatives through a Sobolev estimate and the control in condition \ref{CNLrNice} or \ref{CNLrBetter}. 

We now define the conformal charge density in formal analogy to the corresponding density in $\Reals^{n+1}$, where it is generated by a symmetry of the wave equation. In the present situation, it does not come from a symmetry and does not lead to a conserved quantity. We define the conformal charge to be the integral of the conformal charge density plus the energy density. This slightly strange definition is taken to simplify the statement of certain estimates, for example estimate \eqref{enpKlainermanSobolev}. 
\begin{align}
\ConfChgDens[v,w]
=& (t^2+r_*^2) e[v,w] + 2tr_*p_*[v,w]\hidenumber\\
\ConfChg[v,w]
=& \int_\starman \ConfChgDens[v,w] +e[v,w] \dThree
\hidenumber \end{align}
It has long been known that the conformal charge density can be written as an inherently positive quantity \cite{Strauss}.  
\begin{align}
\label{e12.2}
\ConfChgDens[v,w] 
=& \frac14 (t-r_*)^2(w-v')^2 + \frac14 (t+r_*)^2(w+v')^2 \\
&+ \frac12(t^2+r_*^2)Vv^2 + \frac12 (t^2+r_*^2)V_L(\nabla_{\cpctW}v\cdot \nabla_{\cpctW}v) +\frac12(t^2+r_*^2)f(r_*)F(v^2).
\hidenumber 
\end{align}

In fact, the control on $t^2$ times the angular derivatives is the central element which allows us to prove decay of $L^p$ norms through the Sobolev estimate. The conformal charge can be shown to control the $L^2$ norm, and, interpolating this with the Sobolev estimate provides the following estimate, which is proven in lemma \ref{lpKlainermanSobolev}: for \psianicefn, 
\begin{align}
\int_\starman f(r_*)^3 V_L^{-1} |\psi|^{3(p-1)} \dThree
\leq& C \ConfChg^\frac{3(p-1)}{2} t^{-p+\frac53} .
\hidenumber \end{align}

We now turn to estimating the growth of the conformal charge. Introducing further the Lagrangian and the shear densities
\begin{align}
l[v,w]=& \frac12 (w^2 -v'^2 -Vv^2 -V_L \nabla_{\cpctW}v\cdot\nabla_{\cpctW}v -f(r_*)F(v^2) ) ,\hidenumber\\
T_{**}[v,w]=&v'v' ,\hidenumber\\
\vec{T}_{*\omega}[v,w]=& v'V_L \nabla_{\cpctW} v ,
\hidenumber \end{align}
the relations \eqref{e11.2a}-\eqref{e11.2b} can be rewritten as 
\begin{align}
0=& \dt e - \dr p_* - \nabla_{\cpctW}\cdot \vec{p}_\omega , \hidenumber\\
0=& \dt p_* - \dr l - \dr T_{**} - \nabla_{\cpctW}\cdot\vec{T}_{*\omega} \hidenumber\\
& -\frac12 V'\phi^2 - \frac12 V_L' |\nabla_{\cpctW}\phi|^2 . 
\hidenumber \end{align}
Unless otherwise indicated, we will evaluate all densities on a solution and its $t$ derivative, $(\phi,\dot{\phi})$. From this presentation, and integration by parts, it is easy to show that 
\begin{align}
\dt \ConfChgDens 
=& 2te + 2r_* p_* \hidenumber\\
&+ (t^2+r_*^2) (\dr p_* + \nabla_{\cpctW}\cdot\vec{p}_{\omega}) \hidenumber\\
&+ 2tr_* (\dr l + \dr T_{**} + \nabla_{\cpctW}\cdot T_{*\omega} + \frac12 V'\phi^2 + \frac12 V_L'|\nabla_{\cpctW}\phi|^2 +\frac12 f'(r_*)F(\phi^2) ) , \hidenumber\\
\Dt \ConfChg
=& \int_{\starman} 2t (e-l-T_{**}) + tr_*(V'\phi^2 + V_L'|\nabla_{\cpctW}\phi|^2 + f'F(\phi^2) ) \dThree \hidenumber\\
=& \int_{\starman} t (2V+r_*V')\phi^2 + t (2V_L+r_*V_L')|\nabla_{\cpctW}\phi|^2 + t(2f+r_*f')F(\phi^2) \dThree
\label{e13.4}
\end{align}

To control the growth of the conformal charge, it is sufficient to control a space-time integral which is localized. Introducing $\chia$ to denote a smooth, positive, compactly supported function which dominates $2V+r_*V'$, and $2V_L+r_*V_L'$, for linear solutions, 
\begin{align}
\label{e18.2}
\ConfChg[\phi,\dot\phi](t)
\leq& C\int_0^t \tau \|\chia L \phi\|^2 d\tau + \ConfChg[\phi,\dot\phi](0) 
\end{align}
For nonlinear solutions, the integral $\int_0^t \tau \int_\starman \chi F(|\phi|^2) \dThree d\tau$ must also be controlled where $\chi$ also dominates the nonlinear trapping term, $2f+r_*f'$. 

\section{Three weighted norms controlled by derivatives}
\label{sThreeEstimates}

We present three lemmas in this section. The first is a Sobolev estimate on $\starman$, representing either a warp product manifold of the type we consider or a spatial slice of the exterior region of the Schwarzschild solution. The second is an $L^p$ estimate in terms of the energy and conformal charge, which we need in estimating the nonlinear contributions to the growth of the conformal energy and the growth of the higher angular Sobolev norms. The third is a Hardy or Poincare\'e estimate, which controls a weighted $L^2$ norm by the sum of the $L^2$ norm of the radial derivative and the expectation of the potential. 

We now show that a Sobolev estimate holds both for warp product manifolds of the type we've considered and for the Schwarzschild solution using condition \ref{CNLrNice}. The main feature of this estimate is that it separates the role of the angular derivatives from that of the radial derivatives. This is important, because we will later show that the angular derivatives decay. 

\begin{lemma}
\label{lSobolev}
For \psianicefn, 
\begin{align}
\label{eSobolev}
\int_\starman V_L^2 |\psi|^6 \dThree
\leq C(\int_\starman |\dr\psi|^2 + |(V_L'/V_L) \psi|^2 \dThree) (\int_\starman V_L(|\nabla_{\cpctW}\psi|^2 + |\psi|^2) \dThree)^2 . 
\end{align}
\end{lemma}
\begin{proof}
We reproduce the $H^1\hookrightarrow L^6$ estimate presented in \cite{BlueSoffer3} and following the standard methods for $\Reals^3$ \cite{Evans}. We begin by proving a $W^{1,1}\hookrightarrow L^{\frac32}$ estimate and define
\begin{align}
I_1(\omega) = \int_{\Reals} |\dr \psi_0(r_*,\omega)| dr_*
\hidenumber \end{align}
We use $\psi_0$ to temporarily denote a function with unspecified relation to $\psi$. From the Fundamental Theorem of Calculus, the Cauchy-Schwartz estimate, and the $W^{1,1}(\cpctW)\hookrightarrow L^2(\cpctW)$ Sobolev estimate, for $\psi_0\in\solset$, 
\begin{align}
|\psi_0(r_*,\omega)|^\frac32 \leq& I_1(\omega)^\frac12 |\psi_0(r_*,\omega)| \hidenumber\\
\int_{S^2} |\psi_0(r_*,\omega)|^\frac32 d\omega 
\leq& (\int_{S^2} I_1(\omega) d\omega)^\frac12 (\int_{S^2} |\psi_0(r_*,\omega)|^2 d\omega)^\frac12 \hidenumber\\
\leq& (\int_\starman |\dr\psi_0| \dThree)^\frac12 \int_{S^2} |\nabla_{\cpctW}\psi_0| + |\psi_0| d\omega .
\hidenumber \end{align}
Integrating in $r_*$ against the unspecified weight $h^\frac12$ gives
\begin{align}
\int_\starman h^\frac12 |\psi_0|^\frac32 \dThree \leq (\int_\starman |\dr\psi_0| \dThree)^\frac12 \int_\starman h^\frac12 (|\nabla_{\cpctW}\psi_0| + |\psi_0|) \dThree .
\hidenumber \end{align}
We substitute $\psi_0=h|\psi|^4$ and apply the Leibniz rule and Cauchy-Schwartz, 
\begin{align}
\int_\starman h^2 |\psi|^6\dThree 
\leq& 4(\int_\starman h^2 |\psi|^6\dThree )^\frac34 (\int_\starman |\dr\psi|^2+|h'h^{-1}\psi|^2\dThree)^\frac14 (\int_\starman h|\nabla_{\cpctW}\psi|^2+h|\psi|^2 \dThree)^\frac12 .
\hidenumber \end{align}

We now take $h=V_L$. 
\end{proof}

The following estimate shows that certain $L^p$ norms are controlled by the energy and conformal charge with decaying powers of $t$. The decaying powers of $t$ provide the decay estimates necessary to show that the nonlinear contributions to the growth of the conformal charge and the higher angular Sobolev norms are small. 

\begin{lemma}
\label{lpKlainermanSobolev}
For \psianicefn, under condition \ref{CNLrNice}
\begin{align}
\int_\starman V_L^2 |\psi|^6 \dThree
\leq C \ConfChg[\psi,\dot\psi]^3 t^{-4}
\end{align}
Further, under condition \ref{CNLPotLTheSame}, 
\begin{align}
\label{enpKlainermanSobolev}
\int_\starman f^3 V_L^{-1} |\psi|^{3(p-1)} \dThree
\leq& C \ConfChg^{\frac{3(p-1)}{2}} t^{3(-p+\frac53)} .
\end{align}
\end{lemma}
\begin{proof}
We begin by explaining how to control the inhomogeneous terms, involving no derivatives, on the right of equation \eqref{eSobolev}. This involves controlling the $L^2$ norm by the conformal energy, following \cite{GinibreVelo}, and is similar to the Poincar\'e estimate \eqref{eHardyPoincareSterbenz}. 

Focussing on the time and radial terms in the conformal charge,
\begin{align}
2 \ConfChg[\psi,\dot\psi]
\geq \frac12 \|(t-r)(\dot\psi-\dr\psi)\|^2 + \frac12 \|(t+r_*)(\dot\psi+\dr\psi)\|^2
=\|t\dot\psi+r_*\dr\psi\|^2 +\frac12 \|r_*\dot\psi + t \dr\psi\|^2 .
\hidenumber \end{align}
Adding and subtracting $-h r_*^2(1+r_*^2)^{-1}\psi$ and $h r_*t(1+r_*^2)^{-1}\psi$, and expanding, we have
\begin{align}
\ConfChg
\geq& \frac12\|t\dot\psi+r_*\dr\psi-hr_*^2(1+r_*^2)^{-1}\psi\|^2 +\frac12 \|r_*\dot\psi + t \dr\psi+hr_*t(1+r_*^2)^{-1}\psi\|^2 \hidenumber\\
&-\langle\dr\psi,(t^2-r_*^2)r_*(1+r_*^2)^{-1} h\psi\rangle
-\frac12 \langle h\psi, (t^2+r_*^2)r_*^2(1+r_*^2)^{-2} h\psi\rangle .
\hidenumber \end{align}
Integrating by parts and noting that $r_*$ is differentiable, the last two terms can be written as
\begin{align}
\frac12 \langle \psi, (\frac{r_*^2(t^2+r_*^2)}{(1+r_*^2)^2}(-h-h^2) + \frac{t^2+3r_*^2}{1+r_*^2}h) \psi\rangle .
\hidenumber \end{align}
For $h\in(-1/4,0)$ and $|r_*|>2$, this weight is strictly positive, and using the potential $V$ to control small $r_*$, we have
\begin{align}
\ConfChg[\psi,\dot\psi]
\geq& \frac14 \|(t-r)(\dot\psi-\dr\psi)\|^2 + \frac14 \|(t+r_*)(\dot\psi+\dr\psi)\|^2\hidenumber\\
& + \frac12 \int_\starman (t^2+r_*^2)V \psi^2 \dThree \hidenumber\\
\geq& C \int_\starman  \frac{t^2+r_*^2}{1+r_*^2} \psi^2 \dThree \hidenumber\\
\geq& C\| \psi\|^2 . \label{eL2bound}
\end{align}

By condition \ref{CNLrNice}, we can estimate the $(V_L'/V_L)\psi^2$ term by $\psi^2$ and the $V_L\psi^2$ term by $(1+r_*^2)^{-1}\psi^2$. Since $E\leq\ConfChg$, 
\begin{align}
\label{enKlainermanSobolev}
\int_\starman V_L^2 |\psi|^6 \dThree
\leq& C \ConfChg^3 t^{-4} .
\end{align}
Interpolating this against \eqref{eL2bound}, for $2\leq\sigma\leq6$, 
\begin{align}
\label{enSigmaKlainermanSobolev}
\int_\starman V_L^{\frac{\sigma-2}{2}} \psi^\sigma \dThree
\leq& C \ConfChg^\frac\sigma2 t^{-\sigma+2} .
\end{align}

To conclude, we consider the case $\sigma=3(p-1)$. From condition \ref{CNLPotLTheSame}, $V_L^{\frac{3(p-1)-2}{2}}\geq f^3 V_L^{-1}$, and the desired result holds.
\end{proof}

\begin{remark}
The estimate \eqref{eL2bound} requires that $\psi\rightarrow0$ as $r_*\rightarrow\pm\infty$. These are natural boundary conditions on a warp product, but not entirely natural as $r_*\rightarrow-\infty$ on the Schwarzschild solution. The Sobolev estimate in lemma \ref{lSobolev} still holds for functions which do not decay as $r_*\rightarrow-\infty$, since the radial integration in the proof can be taken from $r_*\rightarrow\infty$, where there is decay. The estimate \eqref{eL2bound} does require decay at $-\infty$; however, by applying the same argument to the region $r_*>0$, and using the contribution from the potential term in the conformal energy, \eqref{eL2bound} can be replaced by $\int_\starman (1-2M/r)^2 \psi^2 \dThree \leq\ConfChg$. This can be interpolated against the Sobolev estimate to control $\int_\starman (1-2M/r)^2 r^{-\sigma+2} \psi^\sigma \dThree$. With $\sigma=3(p-1)$ and $f=(1-2M/r)r^{1-p}$, the integrand is exactly $f^3 V_L^{-1} \psi^{3(p-1)}$. Thus, this lemma would hold for the semi-linear wave equation on the Schwarzschild solution. It is condition \ref{CNLCptTrap}, not \ref{CNLPotLTheSame}, which prevents the main result from holding. 
\end{remark}

\begin{remark}
Under condition \ref{CNLrBetter}, it would be possible to estimate the term $\|r^{-1} \psi\|$ by the Hardy inequality below, so that the weighted $L^\sigma$ norm would be controlled by $C E^{\sigma/6} \ConfChg^{\sigma/3} t^{-\sigma+2}$. 
\end{remark}

The following Hardy or Poincar\'e estimate controls a weighted $L^2$ norm by the $L^2$ norm of the radial derivative plus the expectation value of the potential. We use the contribution from the potential, since our method is to reduce something resembling a one dimensional wave equation, and the radial derivative alone is not sufficient in one dimension. The purpose of this estimate is to show that $\gamma$ is energy bounded. 

\begin{lemma}
\label{lHardyPoincareSterbenz}
For \psianicefn, 
\begin{align}
\label{eHardyPoincareSterbenz}
\|(1+r_*^2)^{-\frac12}\psi\| \leq C E[\psi,\dot\psi] .
\end{align}
\end{lemma}
\begin{proof}
This argument is taken from \cite{BSterbenz}. We start by considering a one dimensional, bounded, $\dot{H}^1$ function, which we will also denote by $\psi$ but with argument $x$, and then integrate the derivative of $(1+|x|)^{-1}\psi^2$ against the weight $\sgn{|x|}$ on a bounded interval.
\begin{align}
\frac{1}{1+|c_1|}(\psi(c_1)^2+\psi(-c_1)^2) - 2\psi(0)^2
=&\int_{-c_1}^{c_1} \sgn(x) \dx(\frac{1}{1+|x|}\psi(x)^2) dx \hidenumber\\
=&\int_{-c_1}^{c_1} \sgn(x)(-\frac{\sgn(x)}{(1+|x|)^2}\psi(x)^2 + \frac{2}{1+|x|}\psi\dx\psi) dx .
\hidenumber \end{align}
Isolating the integral of the weighted term, we have
\begin{align}
\int_{-c_1}^{c_1} \frac{\psi(x)^2}{(1+|x|)^2} dx
=& -\frac{1}{1+|c_1|}(\psi(c_1)^2+\psi(-c_1)^2) + 2\psi(0)^2 + \int_{-c_1}^{c_1}\sgn(x) \frac{2}{1+|x|}\psi\dx\psi dx \hidenumber\\
\leq& 2\psi(0) +(\int_{-c_1}^{c_1}\frac{2}{(1+|x|)^2}\psi^2 dx)^{\frac12} (\int_{-c_1}^{c_1} (\dx\psi)^2 dx)^\frac12 .
\hidenumber \end{align}
Since the weighted norm appears quadratically on the left, but only linearly on the right, taking the limit as $c_1\rightarrow\infty$, 
\begin{align}
\label{e28.3}
\int_{-\infty}^\infty \frac{\psi(x)^2}{(1+|x|)^2} dx
\leq C(\psi(0)^2 + \int_{-\infty}^\infty |\dx \psi|^2 dx) .
\end{align}

On $\starman$, we now generate a family of one dimensional functions by fixing each $\omega$ and $t$ variable and by translating in the radial direction by a quantity $c_0$, that is taking $x=r_0-c_0$. Integrating \eqref{e28.3} over a compact interval of $c_0$ values, 
\begin{align}
\int_{-\infty}^\infty \frac{\psi(r_*)^2}{(1+|r_*|)^2} dr_*
\leq C( \int_{-\infty}^\infty V\psi^2(r_*) + |\dx \psi|^2 dr_* ) .
\hidenumber \end{align}
Integration in the angular variable, and estimating the right hand side by the energy gives the desired result. 
\end{proof}

\section{Heisenberg relation}
\label{sHeisenberg}

We now continue with an approach based on the interpretation of the Schr\"odinger equation as describing the evolution of a quantum system. This allows us to think of certain operators as corresponding to physical observables. We begin by writing the wave equation in terms of an operator, 
\begin{align}
\ddot{\phi} + H\phi = 0
\hidenumber \end{align}
where $H$ is the Hamiltonian defined in terms of the following pieces
\begin{align}
\linH=& H_1+H_2+H_3 ,\hidenumber\\
H_1=& -\drr ,\hidenumber\\
H_2=& V ,\hidenumber\\
H_3=& V_L(-\Delta_{\cpctW}) .
\hidenumber \end{align}
The Hamiltonian, $\linH$, and it pieces, $H_i$, are self-adjoint. From this, we have a Heisenberg relation describing the time evolution of a skew expectation value of an observable. 

\begin{theorem}
If $A(t)$ is a differentiable family of self-adjoint operators with dense common domain, and if \phiasolu, then 
\begin{align}
\label{e9.2}
\Dt(\langle \phi,A\dot\phi\rangle - \langle\dot\phi,A\phi\rangle)
=&\langle \phi,[\linH,A]\phi\rangle + \langle \phi,\dot{A}\dot\phi\rangle - \langle\dot\phi,\dot{A}\phi\rangle +\langle \phi,A f(r_*)F'(\phi^2)\phi\rangle - \langle f(r_*)F'(\phi^2)\phi,A\phi\rangle ,
\end{align}
where $\langle \phi,[\linH,A]\phi\rangle$ is understood to be the quadratic form defined on $D(\linH)\cap D(A)$ by $\langle \linH\phi,A\phi\rangle-\langle A\phi,\linH\phi\rangle$. 
\end{theorem}

This result follows by direct computation. 

Our method for proving local decay and the phase space estimate is to find a propagation observable. A simple example of a \emph{propagation observable} which \emph{majorates} an operator $G$ is a time independent operator $A$ for which 
\begin{align}
[H,A]=G^* G .
\hidenumber \end{align}
If $A$ is \emph{energy bounded}, in the sense that $\|A\phi\|^2\leq E[\phi,\dot\phi]$, then from integrating the Heisenberg relation \eqref{e9.2}, we have a space-time integral estimate on $G$: 
\begin{align}
\int \|G\phi\|^2 dt \leq \int \Dt(\langle\phi,A\dot\phi\rangle-\langle\dot\phi,A\phi\rangle) dt \leq C\|\dot\phi\| \|A\phi\| \leq CE[\phi,\dot\phi] .
\hidenumber \end{align}

The notion of propagation observable is slightly broader in that it allows various lower order terms. 

\begin{definition}
Given a pair of possibly time dependent operators, $A$ and $G$, with $D(A)\cap D(G) \cap D(H)$ dense in $L^2(\starman)$, the operator \emph{$A$ is a propagation observable which majorates $G$} if, for all $\psi\in D(A)\cap D(G) \cap D(H)$, 
\begin{align}
\langle \psi,[\linH,A]\psi\rangle + \langle \psi,\dot{A}\dot\psi\rangle - \langle\dot\psi,\dot{A}\psi\rangle
\geq& \langle G\psi,G \psi\rangle + \langle \psi, R\psi\rangle ,
\hidenumber \end{align}
where $R$ is a sum of \emph{lower order terms} (with respect to $G$), $R_i$, for which either 
\begin{enumerate}
\item $R_i$ is a bounded operator and if \phiasolu, 
\begin{align} 
\int \langle \phi,R_i\phi\rangle dt \leq C
\hidenumber \end{align}
\item or $R_i$ is an operator with $D(R_i)\cap D(G^* G)$ dense and for all $\psi\in D(G^* G)$, 
\begin{align}
|\langle \psi, R \psi\rangle| \leq C\langle \psi, (G^* G)^{1-\epsilon}\psi\rangle
\hidenumber \end{align}
\end{enumerate}
\end{definition}

We make the following additional remarks. The notion of lower order depends on the operators $A$ and $G$. If $A$ is an anti-self-adjoint operator, then the same Heisenberg relation \eqref{e9.2} holds. If $H$ and $A$ take $\solset$ to $\solset$ then the commutator acting on $\solset$ is the operator given by 
\begin{align}
[H,A]=HA-AH .
\hidenumber \end{align}

When making estimates involving operators $B$ and $R$, with $\solset$ dense in their common domain, we will abuse notation by writing
\begin{align}
B \geq R
\hidenumber \end{align}
to mean that for \psianicefn
\begin{align}
\langle \psi,B\psi\rangle \geq \langle \psi, R \psi\rangle .
\hidenumber \end{align}

\section{Local Decay}
\label{sLocalDecay}

The goal of this section is to prove bounds on weighted $L^2$ space-time norms of solutions. We begin by introducing a propagation observable, $\gamma$, show it majorates a weighted quantity, and conclude with a Gronwall's type argument to integrate the Heisenberg-type relation. These estimates are proven using an eigenfunction decomposition analogous to a spherical harmonic decomposition. we refer to this as a harmonic decomposition. We show these estimates have a uniform nature in the harmonic parameter, so that we can recover an estimate for general $u$. In later arguments, the contribution from $H_2$ is controlled by the local decay estimate, and we can use a uniform multiplier. 

In $\Reals^{3+1}$, the radial derivative can be used as a propagation observable to control the time integral of $|u(t,\vec{0})|^2$. This can be thought of as a weighted space-time integral with the $\delta$ function as a weight. Essentially, our propagation observable is a smooth version of the radial derivative, which leads to a smooth weight. 

This follows \cite{LabaSoffer, BlueSoffer1}. Condition \ref{CPtRepulsive} on the effective potential, $\Vtl$, is the main condition used in this section. We recall the effective potential on each spherical harmonic, given by
\begin{align}
\Vtl = H_2 + H_3 = V + \tl^2 V_L .
\hidenumber \end{align}
If necessary, we will use $\Ptl$ to denote the projection on to the harmonic with parameter $l$, but typically, in this section, we will assume that $u$ has support on only one harmonic. In this case, $\phi$ satisfies
\begin{align} 
\ddot{\phi}-\phi''+\Vtl \phi =0 .
\hidenumber \end{align}
From condition \ref{CPtOnePeak}, $\Vtl$ has a single critical point, which is a maximum and which we denote $\al$. 

The Morawetz-type operator is defined on each harmonic in terms of the radial coordinate $r_*$. The definition involves a weight, $g_{l,\sigma}$, which is defined for $\sigma>1$ so that $g_{l,\sigma}$ remains bounded. There is also a dilation factor $b$ which we will later choose sufficiently large, independently of $l$, so that  $\gamma$ majorates $(1+r_*^2)^{-1}$. We will prove this on each harmonic, and then show there is a uniform lower bound on these estimates, to get an estimate uniform in $l$. 

\begin{definition} 
\label{defngamma}
Given $\sigma>1$, $b>0$, and $l\in\Naturals$, the Morawetz-type multiplier $\gamma_{l,\sigma}$ is defined by
\begin{align}
g_{l,\sigma}(r_*)\defin&\int_{0}^{b(r_*-\al)} \frac{1}{(1+|\tau|)^{\sigma}}d\tau\hidenumber\\
\gamma_{l,\sigma} \defin&\frac{1}{2}(g_{l,\sigma}(r_*)\dr+\dr g_{l,\sigma}(r_*))\hidenumber\\
=&g_{l,\sigma}(r_*)\dr +\frac{1}{2}g_{l,\sigma}'(r_*)
\hidenumber \end{align}
In all cases $\sigma$ will not vary so the notation $\gamma=\gamma_{l,\sigma}$ and $g=g_{l,\sigma}$ will be used. When working on more than one harmonic, we also use $\gamma$ to denote 
\begin{align}
\sum_{l} \gamma_{l,\sigma} \Ptl . 
\hidenumber \end{align}
\end{definition}

We note that as $l$ varies, the observable $\gamma_{l,\sigma}$ is simply translated. Since the $\al$ converge, the weights $g_{l,\sigma}$ has a limit in $L^\infty$. 

We now show that the commutators uniformly bound a polynomially decaying term, for an appropriate choice of $b$. This value of $b$ will always be used in $\gamma_{l,\sigma}$ from now on. There are additional positive terms which are also dominated, but we are better able to take advantage of them in the later sections using a slightly modified multiplier. 
\begin{lemma}
\label{lBasicMorawetzCommutatorWasLikeL19andL28}
If $\sigma>1$ and $M>0$, then there is a choice of $b$ for which there is a constant $C$ such that for all $\psi\in\solset$, 
\begin{align}
\label{e5.5}
\langle \psi,[\linH,\gamma]\psi\rangle \geq& C \langle \psi,(1+r_*^2)^{-\frac\sigma2-1}\psi\rangle +C\langle \dr\psi,(1+r_*^2)^{-\frac\sigma2} \dr\psi\rangle    .
\end{align}
\end{lemma}
\begin{proof}
We work on a single spherical harmonic and sum at the end of the argument. The commutator can be computed exactly as 
\begin{align}
[\linH,\gamma]
=& -\frac12 ( \drr g\dr + \drrr g-g\drrr -\dr g\drr ) - (g\dr\Vtl-\Vtl g\dr) \hidenumber\\
=& -2\dr g' \dr -\frac12 g''' - g \Vtl' .
\hidenumber
\end{align}

We will use $x$ to denote $r_*-\al$ 
\begin{align}
\label{engderiva}
g'=& \frac{b}{(1+b|x|)^\sigma} \\
\label{engderivb}
g''=& \frac{-b^2\sigma \sgn(x)}{(1+b|x|)^{\sigma+1}} \\
\label{engderivc}
g'''=& -b^2\sigma 2\delta + \frac{b^3\sigma(\sigma+1)}{(1+b|x|)^{\sigma+2}}
\end{align}

To estimate the contribution from the $\delta$ function, we introduce a smooth, weakly decreasing, compactly support function, $\chi$ which is identically one in a neighborhood of $x=0$ (this is unrelated to any other cut off function used elsewhere in this paper). Using integration by parts, and temporarily treating $\psi$ as a function of $x$, 
\begin{align}
0
=& \int \dr (x\chi \psi^2) dx \hidenumber\\
=&\int \chi \psi^2 dx + \int x\chi'\psi^2dx + \int 2x\chi \psi\dr \psi dx , \hidenumber\\
\label{enForAMi}
\int \chi \psi^2 dx
\leq& \int (|x\chi'|+x^2\chi) \psi^2 dx + \int \chi (\dr \psi)^2 dx ,\\
\psi(0)^2
=& -\int_0^\infty \dr (\chi \psi^2) dx \hidenumber\\
\label{enForAMii}
=& -\int_0^\infty \chi' \psi^2 dx - \int_0^\infty 2\chi \psi\dr \psi dx \\
\leq& \int_0^\infty |\chi'| \psi^2 dx + \int_0^\infty \chi \psi^2 dx + \int_0^\infty \chi (\dr \psi)^2 dx . \hidenumber
\hidenumber \end{align}
Applying a symmetric argument for $(-\infty,0]$ and substituting the result for $\int \chi \psi^2 dx$, we have
\begin{align}
\label{eControlDelta}
2\psi(0)^2 
\leq \int (|\chi'| + |x\chi'|+x^2\chi)\psi^2 dx + \int 2\chi (\dr \psi)^2 dx .
\end{align}

We use a similar integration by parts argument to estimate $(1+b|x|)^{-\sigma-2} \psi$. 
\begin{align}
0=&\int \dr (g''\psi^2) dx \hidenumber\\
=&\int g''' \psi^2 dx + \int 2 g''\psi \dr \psi dx \hidenumber\\
b^3 \int \frac{\sigma(\sigma+1)}{(1+b|x|)^{\sigma+2}}\psi^2 dx
\leq& 2\sigma b^2\psi(0)^2 \hidenumber\\
&+ \left(\int b^3 \frac{\sigma(\sigma+1)}{(1+b|x|)^{\sigma+2}} \psi^2 dx\right)^\frac12 \left(\int\frac{2\sigma}{\sigma+1} b\frac{2(\dr \psi)^2}{(1+b|x|)^\sigma} dx \right)^\frac12 
\label{eResultOfIBP}
\end{align}
Substituting equation \eqref{eControlDelta}, 
\begin{align}
b^3 \int \frac{\sigma(\sigma+1)}{(1+b|x|)^{\sigma+2}}\psi^2 dx
\leq& \sigma b^2 \int (|\chi'| + |x\chi'|+x^2\chi)\psi^2 dx 
+ \sigma b^2 \int 2\chi (\dr \psi)^2 dx \hidenumber\\
&+ \frac12 \int b^3 \frac{\sigma(\sigma+1)}{(1+b|x|)^{\sigma+2}} \psi^2 dx 
+ \int\frac{2\sigma}{\sigma+1} b\frac{2(\dr \psi)^2}{(1+b|x|)^\sigma} dx \hidenumber\\
b^3 \int \frac{\sigma(\sigma+1)}{(1+b|x|)^{\sigma+2}}\psi^2 dx
\leq& 2 \sigma b^2 \int (|\chi'| + |x\chi'|+x^2\chi)\psi^2 dx 
+ 2\sigma b\int (\frac{2}{\sigma+1} \frac{1}{(1+b|x|)^\sigma}+2b\chi) (\dr \psi)^2 dx
\label{enForAMiii}
\end{align} 

The commutator term can now be estimated
\begin{align}
\int 2g'&(\dr \psi)^2 -\frac12 g''' \psi^2 - g\Vtl' \psi^2 dx \hidenumber\\
=& \int \frac{2b}{(1+b|x|)^\sigma} (\dr \psi)^2 dx + \sigma b^2 \psi(0)^2 
-\frac12 \int \frac{b^3 \sigma(\sigma+1)}{(1+b|x|)^{\sigma+2}} \psi^2 dx 
- g\Vtl' \psi^2 dx \hidenumber\\
\geq& \int \frac{2b}{(1+b|x|)^\sigma} (\dr \psi)^2 dx + \sigma b^2 \psi(0)^2 \hidenumber\\
&- \sigma b^2 \int (|\chi'| + |x\chi'|+x^2\chi)\psi^2 dx 
- \sigma b\int (\frac{2}{\sigma+1} \frac{1}{(1+b|x|)^\sigma}+2b\chi) (\dr \psi)^2 dx
- \int g\Vtl' \psi^2 dx \hidenumber\\
\geq& 2b\int (\frac{1}{\sigma+1}\frac{1}{(1+b|x|)^{\sigma+2}} -\sigma b\chi) (\dr \psi)^2 dx 
+ \int (-\sigma b^2 (|\chi'| + |x\chi'|+x^2\chi) - g\Vtl') \psi^2 dx 
+ \sigma b^2 \psi(0)^2 
\label{enForAMiv}
\end{align}
The term $\sigma b \psi(0)^2$ is positive and will be ignored. The term $(\sigma+1)^{-1}(1+b|x|)^{-\sigma-2} -\sigma b \chi$ is independent of $l$, so we can choose $b$ sufficiently small so that this is bounded below by $c(1+b|x|)^{-\sigma-2}$. Since $\Vtl' = V' + V_L'\tl^2$, and the $\al$ converge, the value of $\Vtl''$ is uniformly bounded from below at $\al$, and $b$ can be chosen sufficiently small so that the quadratically vanishing quantity $-\sigma b^2 (|\chi'| + |x\chi'|+x^2\chi) - g\Vtl') $ is bounded below by $\sigma b^2 (|\chi'| + |x\chi'|+x^2\chi) $, uniformly in $l$. Thus there is a uniform constants, $c_1$, for which, 
\begin{align}
\int 2g'&(\dr \psi)^2 -\frac12 g''' \psi^2 - g\Vtl' \psi^2 dx \hidenumber\\
\geq& c_1 \int \frac{1}{(1+b|x|)^\sigma}(\dr \psi)^2 + \sigma b^2 (|\chi'| + |x\chi'|+x^2\chi) \psi^2 dx \hidenumber
\end{align}
Applying the estimate \eqref{enForAMiii}, there is a uniform constant, $c_2$, such that
\begin{align}
\int 2g'&(\dr \psi)^2 -\frac12 g''' \psi^2 - g\Vtl' \psi^2 dx \hidenumber\\
\geq& c_2 \int \frac{1}{(1+b|x|)^{\sigma+2}} \psi^2 dx  + c_2 \int \frac{1}{(1+b|x|)^{\sigma}} (\dr\psi)^2 dx.
\label{e5.4}
\end{align}
Since the $\al$ converge, there is a uniform equivalence between $(1+b|x|)^{-\sigma-2}$ and $(1+r_*^2)^{-(\sigma/2)-1}$. Integrating over the angular variables gives
\begin{align}
\langle \psi, [\linH,\gamma]\psi\rangle
\geq& C \langle \psi, \frac{1}{(1+r_*^2)^{\frac\sigma2+1}} \psi\rangle 
+ C \langle \dr\psi, \frac{1}{(1+r_*^2)^{\frac\sigma2}} \dr\psi\rangle .
\hidenumber \end{align}
\end{proof}

For the local decay result, equation \eqref{e5.5} is sufficient, but in the phase space analysis, we will need more contributions from $-g V_L'$. We now prove the local decay estimate. Since $\gamma$ is self-adjoint and preserves the class of real valued functions 
\begin{align}
\label{e6.1}
\langle \dot\psi,\gamma\psi\rangle-\langle\psi,\gamma\dot\psi\rangle =& 2\langle\dot\psi,\gamma\psi\rangle , \\
|\langle \dot\psi,\gamma\psi\rangle-\langle\psi,\gamma\dot\psi\rangle| \leq& 2\|\dot\psi\| \|\gamma\psi\|\hidenumber\\
\leq& C\|\dot\psi\| \|\dr\psi\| + C\|\dot\psi\| \|(1+r_*^2)^{-\sigma/2}\psi\| .
\hidenumber \end{align}

From lemma \ref{lHardyPoincareSterbenz}, the energy controls the weighted term arising from $\gamma$, 
\begin{align}
 \|(1+r_*^2)^{-\sigma/2}\psi\| 
\leq \|(1+(r_*)^2)^{-\frac12}\psi\| 
\leq C E[\psi,\dot\psi] .
\label{eGammaEnergyBound}
\end{align}

This proves that $\gamma$ is energy bounded. Taking $\sigma=2>1$ for simplicity, integrating the Heisenberg relation, and using \eqref{e5.5} and \eqref{eGammaEnergyBound}, we have that for a solution to the linear wave equation,
\begin{align}
\int_0^t \langle \phi,[H,\gamma]\phi\rangle d\tau = 2\langle\phi,\gamma\phi\rangle |_0^t \label{e7.1}\\
\int_0^t \| (1+r_*^2)^{-1} \phi\|^2 d\tau \leq& CE[\phi,\dot\phi] .
\hidenumber 
\end{align}

The nonlinear contribution from the Heisenberg formula is
\begin{align}
|- \langle \phi,\gamma fF'\phi\rangle+ \langle fF'\phi,\gamma\phi\rangle|
\leq E[\phi,\dot\phi]^\frac12 \| fF'\phi\|  \hidenumber 
\end{align}

Using this result and estimate \eqref{e5.5}, we can conclude
\begin{align}
C E[\phi,\dot\phi]\geq \int_0^t  
\| (1+r_*^2)^{-1} \phi\|^2 
+ \|(1+r_*^2)^{-\frac12}\dr \phi\|^2 
- E[\phi,\dot\phi]^\frac12 \| fF'\phi\|
d\tau .
\label{eLocDecWithOutGVL}
\end{align}
This proves equation \eqref{eLocDecInIntro} in the linear case. 

We now rerun the argument with
\begin{align}
g=\int_0^{b(r_*-(\alpha_\infty)_*)} \frac{1}{(1+|\tau|)^\sigma} d\tau ,
\label{eOneg}
\end{align}
In the term $-g\Vtl'=-gV'+g\tl^2V_L'$, we treat $-gV'$ as an error term to be controlled and $-g\tl^2V_L'$ as a good term. The term$-gV'$ is negative between $(\alpha_0)_*$ and $(\alpha_\infty)_*$, which is a compact region, so $-gV' \psi^2$ can be controlled by $(1+r_*^2)^{-1} \psi^2$ in the previous estimate. We can now use $-g \tl^2 V_L'$ as in the previous estimate to get positivity and to get an extra $-g \tl^2 V_L'$ term, so that we have \eqref{eLocDecWithOutGVL} with an extra term. 
\begin{align}
C E[\phi,\dot\phi]\geq \int_0^t (& 
\| (1+r_*^2)^{-1} \phi\|^2 
+ \|(1+r_*^2)^{-\frac12}\dr \phi\|^2 \hidenumber\\&
+ \|(- g V_L)^\frac12 \nabla_{\cpctW}\phi\|^2 
- E[\phi,\dot\phi]^\frac12 \| fF'\phi\| )
d\tau .
\hidenumber
\end{align}

\section{Angular refinements}
\label{sAngularRefinements}

We use two methods to gain additional angular regularity: angular modulation and phase space localization. Angular modulation involves rescaling the radial variable using the angular Laplacian operator. For phase space localization, we localize the observables in regions where the radial derivative is comparable to powers of the angular derivative. The angular modulation alone allows us to regain control over $\frac34$ derivatives in $\Lebesgue^2$. From the local decay result, we already control a full derivative in the radial direction. The phase space localization allows is to work in ``bands'' where powers of the angular derivative are related to the radial derivative. This allows us to translate the radial derivative control in to control of an angular derivative with the loss of $\delta$ derivatives from the width of the ``bands''. There is also a $\epsilon$ loss from the shape of the localization function. We combine these as a $\varepsilon$ loss of angular derivatives. 


To compute commutators involving the phase space localization, we use a commutator expansion lemma. To use this we need to introduce a notion of lower order, which we will take to mean local terms involving fewer powers of the angular derivative. 

\begin{definition}
The operator $L$ is defined by the spectral theorem as 
\begin{equation}
L = \sqrt{1-\wlap}
\end{equation}

The phase space variables are defined by 
\begin{align}
\qm =& L^m r_* \hidenumber\\
\pn =& -i L^{n-1} \dr
\hidenumber \end{align}
\end{definition}

\begin{definition}
The phrase ``lower order terms'' will refer to terms of the form 
\begin{align}
\| L^{\frac12+n} \chi u \|^2 ,
\hidenumber \end{align}
and terms which are controlled by the local decay estimate such as
\begin{align}
\| \chi \psi\|^2, && \| \chi \dr \psi \|^2 , && \| \nabla_{\cpctW}\tilde\chi \psi\|^2, 
\hidenumber \end{align}
where $\chi$ is a compactly supported function, and $\tilde\chi$ is a compactly supported function which vanishes at least linearly at $r_*=0$. These are denoted 
\begin{align}
\ErrorTermsnPsi
\hidenumber \end{align}
\end{definition}

The estimate for the phase space observable involves localization in the phase space variables. For this purpose, the following notation is used when the localizing function is a characteristic function.  

\begin{definition}
For $\delta>0$, the sharp, near localization and sharp, interval localization for $\pn$ are 
\begin{align}
\pnn{\xi} =& \chi([0,1],\xi) \hidenumber\\
\pniWITHl{\xi} =& \chi([l^{-\delta},1],\xi) .
\hidenumber \end{align}
From these, $\pnn{\pn}$ and $\pni{\pn}$ can be defined by the spectral theorem. 
\end{definition}

We now define the phase space observable, $\Gammahat$, and its pieces, $\Gammanm$. To do this, we will first need to introduce a more regular phase space localization function. Recall that $\chia$ is defined to be a positive, smooth, compactly supported function which dominates the trapping terms. Essentially, it localizes near $r=\alpha$. 

\begin{definition}
For fixed values of $\delta>0$ and $\epsilon>0$, the following definitions are made. For simplicity, we will take $\delta$ to be the inverse of a positive, even integer. 

The angularly modulated Morawetz multiplier is defined in terms of $g$ from \eqref{eOneg} by 
\begin{align}
\gm=& g(L^m r_*) = g(\qm) ,\hidenumber\\
\gammam=&\frac12(\gm\dr +\dr\gm) .
\hidenumber 
\end{align}

The first smooth, near phase space localization functions is defined by 
\begin{align}
\pnna{x} =& (1+x^2)^{-\frac{1-\epsilon}{4}} .
\hidenumber \end{align}
The partial phase space observables are defined for $0\leq n \leq m \leq \frac12$ by 
\begin{align}
\Gammanm 
=& \chia \pnna{\pn} L^{n-\epsilon} \gamma_{L^m} \pnna{\pn} \chia \hidenumber\\
=& \frac12 \chia \pnna{\pn} i \left(\gm\pn + \pn \gm\right) \pnna{\pn} \chia
\hidenumber \end{align}

The phase space observable is defined for $\delta=\epsilon$ by 
\begin{align}
\Gammahat=& \Gamma_{\frac12-2\delta,\frac12} + \sum_{j=0}^{\frac{1}{2\delta}-2}\Gamma_{j\delta,j\delta}
\hidenumber \end{align}
\end{definition} 
 
With the local decay observable, the phase space observable majorates $L^{1-\varepsilon}$. Alone the commutator of the phase space observable dominates $1-\varepsilon$ angular derivatives localized for $|-i\dr|\leq L$, up to lower order terms. The remaining regions of phase space and the lower order terms can be controlled by the local decay observable. For our later arguments, the estimate on the phase space observable alone will be necessary. 

\begin{lemma}
\label{lGammahatMajoratesL1minusepsilon}
For \psianicefn, $\delta>0$, there is a $\CGamma>0$ such that
\begin{align}
\langle \psi, [H,\Gammahat] \psi\rangle
\geq& C \| L^{1-\frac{3\delta}{2}} \pnn{-iL^{-1}\dr} \chi \psi \|^2 + \ErrorTermsMinusTwoDeltaPsi \hidenumber\\
\langle \psi, [H,\Gammahat + \CGamma \gamma] \psi\rangle
\geq& C (\| L^{1-\frac{3\delta}{2}} \chi \psi \|^2 + \|(1+r_*^2)^{-1}\psi\|^2)
\label{ePre4.14}
\end{align}
\end{lemma}
\begin{proof}
Each partial phase space observable majorates a power of $L$ in a region of phase space. The partial phase space observables appearing in the phase space observable are chosen so that, together, they cover the entire phase space, up to lower order terms. Since the lower order terms are controlled by the local decay observable, it is sufficient to add a sufficiently large multiple of the local decay observable to get a strictly positive commutator. 

The estimates on the partial phase space observables are reproduced from \cite{BlueSoffer3} in appendix \ref{saPhaseSpaceEst}. Lemmas \ref{lEstimatingThePhaseSpaceCommutatorBeneathStrip} and \ref{lEstimatingThePhaseSpaceCommutatorInStrip} say that the partial phase space observables with $m=n$ dominate $L^{1-\delta-\frac\epsilon2}$ in strips of width $L^\delta$ in $\pn$, and that the partial phase space observables with $n=\frac12$ dominate $L^{\frac34+\frac{n-\epsilon}{2}}$ for $|\pn| < 1$. 
\begin{align}
\langle \psi,[\linH,\Gammanhalf]\psi\rangle\geq& C\|L^{\frac{3/2+n-\epsilon}{2}}\pnn{\pn}\chia \psi\|^2-\ErrorTermsnPsi \hidenumber\\
\langle \psi,[\linH,\Gammann]\psi\rangle\geq& \|L^{1-\delta-\frac\epsilon2}\pni{\pn}\chia \psi\|^2 -\ErrorTermsnPsi
\hidenumber \end{align}

Summing over the partial phase space observables, $\Gammanm$, in the phase space observable, $\Gammahat$, gives an estimate in many regions of phase space. The highest powers of $L$ appearing in the lower order terms are $\frac{1+2n}{2}$ for $n=\frac12-2\delta$. 
\begin{align}
\langle \psi,[H,\Gammahat]\psi\rangle
\geq& C (\| L^{1-\frac{3\delta}{2}} \pnn{\phalfminustwodelta} \chia \psi \|^2
+ \sum_{j=0}^{\frac{1}{2\delta}-2} \| L^{1-\frac{3\delta}{2}} \pni{\pjdelta} \chia \psi \|^2)
 + \ErrorTermsMinusTwoDeltaPsi
\hidenumber \end{align} 

The regions $|\phalfminustwodelta|\leq1$ and $L^{-\delta}\leq|\pjdelta|\leq 1$ correspond to $|-i\dr|\leq L^{\frac12 + 2\delta}=L^{1-(\frac12-2\delta)}$ and $L^{1-(j+1)\delta}\leq |-i\dr|\leq L^{1-j\delta}$. Since $j$ varies from $0$ to $\frac{1}{2\delta}-2$, the regions covered, taken together, cover all of the region $|-i\dr|\leq L$. 
\begin{align}
\langle \psi,[H,\Gammahat]\psi\rangle
\geq& C \| L^{1-\frac{3\delta}{2}} \pnn{-iL^{-1}\dr} \chia \psi \|^2 +\ErrorTermsMinusTwoDeltaPsi
\hidenumber \end{align}

The region $|-i\dr|\geq L$ is covered by the local decay observable. 
\begin{align}
\langle\psi,[\linH,\gamma]\psi\rangle
\geq& \| \chia \dr \psi\|^2 + \|\chia' \psi\|^2 \hidenumber\\
\geq& \| \dr \chia \psi\|^2 \hidenumber\\
\geq& \| L (1-\pnn{-iL^{-1}\dr}) \chia\psi\|^2 .
\hidenumber \end{align}
Putting the estimates in the different regions of phase space together, we have that
\begin{align}
\langle \psi,[H,\Gammahat + C\gamma]\psi\rangle
\geq& C \|L^{1-\frac{3\delta}{2}} \chia \psi \|^2 
+ \ErrorTermsMinusTwoDeltaPsi .
\hidenumber \end{align} 
Since $1-\frac{3\delta}{2} > 1-2\delta$, the first term dominates the lower order terms for sufficiently large $L$. Since the local decay observable majorates $\|(1+r_*^2)^{-1}\psi\|^2 \geq C\|\chi \psi\|^2$ with no remainder, by interpolation, there's a constant $\CGamma$ for which 
\begin{align}
\langle \psi,[H,\Gammahat + \CGamma\gamma]\psi\rangle
\geq& C (\|L^{1-\frac{3\delta}{2}} \chia \psi \|^2 
+ \|(1+r_*^2)^{-1}\psi \|^2) .
\hidenumber \end{align}
\end{proof}

This result can now be integrated, using the Heisenberg relation to conclude the $L^2$ space-time integrability result for $1-\varepsilon$ angular derivatives. 

\begin{proposition}
\label{P15.1}
If \phialinsolu{} and $\varepsilon>0$
\begin{align}
C\int_0^\infty \|L^{1-\varepsilon}\chia\phi\|^2 d\tau
\leq& E[\phi,\dot\phi]
\hidenumber \end{align}
\end{proposition}
\begin{proof}
As in the proof of equation \eqref{e7.1}, we apply the Heisenberg relation and integrate the commutator estimate \eqref{ePre4.14}. To do this, we need to show $\Gammahat+\CGamma\gamma$ is energy bounded. From \eqref{e6.1}, we already have a bound on the term coming from $\gamma$ at the initial and final time. Since $\Gammahat$ is a finite sum of $\Gammanm$ terms, it is sufficient to estimate each $\Gammanm$ term individually. From lemma \ref{lGammanmIsEBndd}, for any function $\psi\in\solset$, 
\begin{align}
\|\Gammanm\psi\|^2 \leq C_{n,m} E[\psi,\dot\psi]
\hidenumber \end{align}
We can now integrate the Heisenberg relation with $A=\Gammahat+\CGamma\gamma$. 
\begin{align}
\label{e16.2}
\Dt(\langle\dot\phi,A\phi\rangle-\langle\phi,A\dot\phi\rangle)
=&\langle\phi,[\linH,A]\phi\rangle \\
2\langle\dot\phi,(\Gammahat+\CGamma\gamma)\phi\rangle|_0^t
\geq&\int_0^t \|L^{1-\frac32 \delta}\chia\phi\|^2 + \|(1+r_*^2)^{-1}\phi\|^2 d\tau . 
\hidenumber \end{align}
Estimating the right hand side by the energy, using the Hardy type estimate \eqref{lHardyPoincareSterbenz}, and taking the limit as $t\rightarrow\infty$. 
\begin{align} 
C E[\phi,\dot\phi]
\geq& \int_0^\infty \|L^{1-\frac32 \delta}\chia\phi\|^2 + \|(1+r_*^2)^{-1}\phi\|^2 d\tau . \hidenumber 
\end{align}
We take $\varepsilon=(3/2)\delta$. 
\end{proof}

\section{Temporal Refinements}
\label{sTemporalRefinements}

We now apply the method from \cite{BSterbenz} to gain additional control in time. Instead of introducing an energy bounded observable, we introduce an observable which is controlled by the conformal energy, and use a boot strap argument to control the growth of the conformal energy. 

We begin by introducing a localization deep inside the light cone and another concept of lower order based on this localization. 

\begin{definition}
Let $\chitfn:\Reals\rightarrow[0,1]$ be a smooth, even function which is identically one in $[-1,1]$, identically zero outside $[-2,2]$ and decreases smoothly in $[1,2]$. If no argument is given, it is assumed that
\begin{align}
\label{e24.1}
\chit= \chitfn(\frac{10r_*}{1+t}) . 
\end{align}
The constant $10$ in this definition is not particularly relevant. With this, or other choices of constant greater than one, $\chit$ will be called the \emph{light cone localization}. 

The phrase ``\emph{lower order in the light cone}'' will refer to terms of the form, or bounded by, 
\begin{align}
E[\chit\psi,\chit\dot\psi] ,
\hidenumber \end{align}
and will be denoted
\begin{align}
\OSterbenz .
\hidenumber \end{align}
\end{definition}

Note that, because we allow different constants to appear in the definition of $\chit$, and 
\begin{align}
\dr(\chitfn(\frac{10r_*}{1+t}))
=\frac{10}{1+t}\chitfn'(\frac{10r_*}{1+t})
\leq \frac{C}{1+t}\chitfn(\frac{5r_*}{1+t}),
\hidenumber \end{align}
we write $\dr \chit = (1+t)^{-1}\chit$. As long as only finitely many derivatives of $\chit$ are taken, we can control $\chit$ and its derivatives by $\chit$ taken with a different constant in the argument. Each time a derivative is taken, because the support of $\chitfn$ is in $|x|\leq2$, the constant in the argument decrease by a factor of $2$. By starting with a sufficiently large constant in equation \eqref{e24.1}, we can ensure that, even after dominating derivatives, the constant in $\chit$ will be larger than $2$, and the support of $\chit$ will be strictly inside the light cone. In particular, we will never take more than two derivatives, so the choice of constant $10$ in \eqref{e24.1} is sufficient. Were we to require control on more derivatives of $\chit$, we would simply take a larger constant in \eqref{e24.1}. 

By the Poincar\'e type estimate in lemma \ref{lHardyPoincareSterbenz}, we can also control weighted norms of $\phi$ with light cone localization. 
\begin{align}
\int \frac{\psi^2}{1+r_*^2} \dThree\leq& C E[\psi,\dot\psi]\hidenumber\\
\int \frac{1}{1+r_*^2} (\chit\psi)^2 \dThree \leq& CE[\chit\psi,\chit\dot\psi]=\OSterbenz .\hidenumber 
\end{align}
This also means that the order of $\chit$ localization and derivatives does not matter. 
\begin{align}
\dr (\chit\psi)
=&\chitfn(\frac{10r_*}{1+t})\dr\psi+\frac{10}{1+t}\chitfn'(\frac{10r_*}{1+t})\psi \hidenumber\\
\leq&\chitfn(\frac{10r_*}{1+t})\dr\psi+\frac{C}{1+t}\chitfn(\frac{5r_*}{1+t})\psi \hidenumber\\
\leq&\chitfn(\frac{10r_*}{1+t})\dr\psi+\frac{C}{(1+r_*^2)^{\frac12}}\chitfn(\frac{5r_*}{1+t})\psi \hidenumber\\
\| \dr (\chit\psi) \|^2 =&\OSterbenz . \hidenumber 
\end{align}

Using this light cone localization, we introduce the temporal phase space observable. The first part is $(1+t)$ times the phase space observable, and the second is $(1+t)$ times the local decay observable localized inside the light cone. There is no need to localize the phase space observable inside the light cone since it already contains the localization $\chia$, which is localized near the photon sphere.

\begin{definition}
The temporal phase space observable is 
\begin{align}
\Gammacrown
=& (1+t) \Gammahat + C_\Gamma (1+t)\chit(\frac{10 r_*}{1+t}) \gamma \chit(\frac{10 r_*}{1+t}) \hidenumber 
\end{align}
\end{definition}

Note that $\Gammacrown$ is anti self adjoint since $1+t$ and $\chit$ are bounded functions at each time, and $\Gammahat$ and $\gamma$ are anti-self-adjoint. 

The temporal phase space observable majorates $1-\varepsilon$ angular derivatives and $\frac12$ factors of $t$; that is, the commutator dominates $2-2\varepsilon$ angular derivatives and one factors of $t$ in expectation value. To show this, we must also include the contribution from the time derivative. 

\begin{lemma}
\label{lGammacrownMajoratesLOneMinusEpsilon}
For \psianicefn, and $\varepsilon>0$, 
\begin{align}
\langle \psi, [H,\Gammacrown]\psi\rangle + \langle \psi, \dot{\Gammacrown} \dot{\psi}\rangle - \langle \dot{\psi},\dot{\Gammacrown}\psi\rangle
\geq& (1+t)C( \| L^{1-\varepsilon} \chi \psi \|^2 + \| (1+r_*^2)^{-1}\chit \psi\|^2) 
+ \OSterbenz
\hidenumber \end{align}
\end{lemma}
\begin{proof}
This lemma says that $\Gammacrown$ majorates $t^\frac12 L^{1-\epsilon}$. Since $\Gammacrown$ is essentially $t\Gammahat$, this should be expected. In essence, the left hand side is $t$ times the commutator with $\Gammahat$, so that the result follows from proposition \ref{lGammahatMajoratesL1minusepsilon}. There are corrections from the commutator with the light cone localization in $\Gammacrown$ and from the time derivative of $\Gammacrown$. These are shown to be lower order inside the light cone. 

The commutator can be broken into four pieces using the linearity of the commutator and the Leibniz rule for commutators. 
\begin{align}
\label{eHGammacrowncommutator}
[H,\Gammacrown]
=&(1+t)[H,\Gammahat] + \CGamma(1+t)\chit[H,\gamma]\chit + \CGamma(1+t)[H,\chit]\gamma \chit + \CGamma(1+t)\chit\gamma[H,\chit]
\end{align}
We now show that the first two terms provide the bound that we want, that the three remaining terms are lower order in the light cone, and that the terms involving $\dot\Gammacrown$ are also lower order in the light cone. Combining these three steps proves the desired result. 

\newstepsequence
\newstep{The dominant terms in the commutator}
The contribution from the first two terms can be found using proposition \ref{lGammahatMajoratesL1minusepsilon}. 
\begin{align}
\langle \psi, (1+t)[H,\Gammahat] + \CGamma(1+t)\chit[H,\gamma]\chit \psi \rangle
\geq (1+t)\Big(& C\| L^{1-\frac32 \delta} \pnn{-iL^{n-1}\dr} \chia u\|^2 + \ErrorTermsMinusTwoDeltaPsi \hidenumber \\
&+ \CGamma \langle\psi,\chit[H,\gamma]\chit\psi\rangle \Big) \hidenumber 
\end{align}
Since $(1+r_*^2)^{-1} \chit > C\chia$, the interpolation in proposition \ref{lGammahatMajoratesL1minusepsilon} still holds, and 
\begin{align}
\langle \psi, (1+t)[H,\Gammahat] + \CGamma(1+t)\chit[H,\gamma]\chit\chit \psi \rangle
\geq& (1+t)C \left( \| L^{1-\varepsilon} \chia \psi \|^2 + \| (1+r_*^2)^{-1}\chit(\frac{10 r_*}{1+t}) \psi\|^2 \right) \hidenumber\\
\geq& (1+t)C \| L^{1-\varepsilon} \chia \psi \|^2 .
\hidenumber \end{align}

\newstep{Remaining terms in the commutator}
We now consider the last two terms in equation \eqref{eHGammacrowncommutator}. The fourth term is, up to signs, the adjoint of the third. When the derivative in the commutator lands on $\chit$, a negative power of $(1+t)$ is introduced, so that the commutator is lower order in time. It is also localized deep in the light cone. 
\begin{align}
(1+t)[H,\chit]
=& (1+t)\left( -2\dr(\dr \chit(\frac{10 r_*}{1+t})) + (\drr \chit(\frac{10 r_*}{1+t})) \right) \hidenumber\\
=& \left( -20 \dr \chit'(\frac{10 r_*}{1+t}) - \frac{100}{1+t} \chit''\right) \hidenumber\\
\gamma\chit=&\chit\gamma + \frac{10 g }{1+t} \chit' 
\hidenumber \end{align}
Since $\gamma\chit$ and the adjoint of $(1+t)[H,\chit]$ both look like the sum of terms involving one or zero radial derivatives times bounded functions and light cone localized functions, the remaining commutator terms are lower order inside the light cone. 
\begin{align}
\langle \psi, (1+t)[H,\chit]\gamma\chit \psi\rangle + \langle \psi, (1+t)\chit \gamma [H,\chit] \psi\rangle
=& \OSterbenz .
\hidenumber \end{align} 

\newstep{Controlling the $\dot{\Gammacrown}$ terms}
\label{sGammacrowndot}
Since $\Gammacrown$ is anti-self-adjoint, it is sufficient to consider 
\begin{align}
\langle \dot\psi, \dot{\Gammacrown}\psi\rangle .
\hidenumber \end{align}

We first compute $\dot{\Gammacrown}$. 
\begin{align}
\dot{\Gammacrown} = \Gammahat + \CGamma \chit\gamma\chit - \CGamma(1+t)\frac{10r_*}{(1+t)^2}(\chit'\gamma\chit+\chit\gamma\chit') .
\hidenumber \end{align}
We will deal with each piece separately. 

Since $\chia$ is localized inside a fixed compact region, and $\chit$ is identically one in a larger region for all $t\geq1$, $\chit$ localization can be added to $\Gammahat$ for free. 
\begin{align}
\Gammahat=\chit\Gammahat
\hidenumber \end{align}
Therefore, 
\begin{align}
\langle \dot\psi, \Gammahat \psi\rangle
=& \langle \chit \dot\psi, \Gammahat \psi\rangle\hidenumber\\
\leq& \|\chit \dot\psi\| \| \Gammahat \psi\|
\hidenumber \end{align}
From lemma \ref{lGammanmIsEBndd}, $\| \Gammahat\psi\|^2$ is controlled by the energy of $\chia\psi$, so 
\begin{align}
\langle \dot\psi, \Gammahat \psi\rangle
\leq& \|\chit \dot\psi\| E[\chia \psi]^\frac12 \hidenumber\\
\leq& \|\chit \dot\psi\|^2 + E[\chit\psi]\hidenumber\\
\leq& \OSterbenz
\hidenumber \end{align}

Clearly, the $\chit\gamma\chit$ terms are lower order inside the light cone. 
\begin{align}
\langle\dot\psi,\chit\gamma\chit\psi\rangle
\leq& \|\chit\dot\psi\| \|\gamma\chit\psi\| = \OSterbenz
\hidenumber \end{align}

The remaining terms, from when the time derivative falls on $\chit$, are dealt with similarly. The derivative gives an extra factor of $r_*/(1+t)$. This is bounded inside the light cone. 
\begin{align}
\langle\dot\psi,(1+t)\frac{10r_*}{(1+t)^2}(\chit'\gamma\chit+\chit\gamma\chit')\psi\rangle
\leq& |\langle \chit'\dot\psi, 10  \frac{|r_*|}{1+t}\gamma\chit\psi\rangle| + |\langle \chit\psi, 10\frac{|r_*|}{1+t} \gamma\chit'\psi\rangle |
\leq \OSterbenz
\hidenumber \end{align}
\end{proof}

In our proof, we bound terms both by the energy and the conformal charge. Energy boundedness was central for proving equation \eqref{e7.1}, and we require an analogue combining $E$ and $\ConfChg$. For \psianicefn, $Q(t)$ depending on $t$ and $\psi$ is \emph{conformally bounded} if 
\begin{align}
Q(t)\leq \min\{(1+t)E[\psi,\dot\psi],(1+t)^{-1}\ConfChg[\psi,\dot\psi] \} .\hidenumber 
\end{align}
We say an operator $A$ is conformally bounded if $\|A\psi\|^2$ is conformally bounded. 

We now turn to proving the main estimate, that the conformal charge of $\phi$ and the energy of $L^\varepsilon\phi$ remain bounded for the linear wave equation or for the nonlinear wave equation with small initial data. In the linear case, we can apply $L^\varepsilon$ to the wave equation to show that the energy of $L^\varepsilon\phi$ is conserved, and then apply the angular refinement argument to control $L\phi$ in $\Lebesgue^2$. This controls the growth of the conformal charge. In the nonlinear case, we can still use the energy of $L^\varepsilon\phi$ to control the growth of the conformal charge. The conformal charge can then be used to control weighted $\Lebesgue^q$ norms, which in turn control the growth of the energy of $L^\varepsilon\phi$. Through a bootstrap, for small data, the energy of $L^\varepsilon\phi$ and the conformal charge control each other. 

\begin{proposition}
If $\varepsilon>0$ and \phialinsolu, 
\begin{align}
\ConfChg[\phi(t),\dot\phi(t)]
\leq& C \ConfChg[\phi(0),\dot\phi(0)] + C E[L^\varepsilon\phi(0),L^\varepsilon\dot\phi(0)] 
\hidenumber \end{align}

For $8/3<p<3$ and sufficiently small initial data, in the norms appearing on the right hand side, the same estimate is true if \phiasolu. 
\end{proposition}
\begin{proof}
If it weren't for the loss of $\varepsilon$ angular derivatives on the commutator estimate of lemma \ref{lGammacrownMajoratesLOneMinusEpsilon}, we'd expect to be able to close the estimate on the growth of the conformal charge using conformal boundedness of the operator $\Gammacrown$. Instead, we will have to resort to a more complicated notion of $L^\varepsilon$-conformal boundedness. 

Our method follows from the following six estimates, in terms of a constant $N$ to be chosen later, 
\begin{align}
\label{n2.1}
\ELeighteen(t)
\leq& \ELeighteen(0) &&\text{Persistence of regularity} \\
&+ \int_0^t  \ELeighteen^\frac12 \| L^{18\varepsilon} f(r_*) F'(|\phi|^2) \phi \|  d\tau
&& \notag\\
\label{n2.2}
\ConfChgphi(t)
\leq& \ConfChgphi(0) + C\int_0^t \tau \|L\chia\phi\|^2 d\tau 
&&\text{Conformal estimate}\\
&+ C\int_0^t \tau \int_\starman (2f+r_*f')F(|\phi|^2) \dThree d\tau 
&&\notag\\
\label{n2.3}
\int _0^t \tau \|L\chia\phi\|^2 d\tau
\leq& C N^2 \supt \ELeighteen(\tau) 
&&\text{Temporal refinement}\\
& + C N^{-\frac14}\supt\ConfChg(\tau) 
+C \int_0^t |\langle fF'\phi,L^{2\varepsilon}\Gammacrown \phi\rangle | d\tau
&&\notag\\
\label{enNLPersistanceOfRegularity}
\int_0^t \ELeighteen^\frac12 &\| L^{18\varepsilon} f(r_*) F'(|\phi|^2) \phi \|  d\tau
&&\text{NL persistence error}\\
\leq& C \supt \ELeighteen(\tau) \supt \ConfChgphi(\tau)^\frac{p-1}{2}  
&& \notag\\
\label{enNLConformal}
\int_0^t \tau \int_\starman (2f+r_*f')&F(|\phi|^2) \dThree d\tau 
&&\text{NL conformal error}\\
\leq& C N^{\frac12} E + C N^{\frac{-p+1}{4}} \supt\ConfChg(\tau)^{\frac{p+1}{2}} 
&&\notag\\
\label{enNLMorawetz}
\int_0^t \langle fF'\phi,L^{2\varepsilon}\Gammacrown \phi\rangle d\tau
\leq& C\supt \ELeighteen(\tau)^\frac12 \supt \ConfChgphi(\tau)^\frac{p}{2}
&&\text{NL Morawetz error}
\end{align}
The conformal estimate \eqref{n2.2} is already known from \eqref{e13.4}. 

In the linear case, the persistence of regularity is trivial since $L^{18\varepsilon}$ commutes with the d'Alembertian, so the $\ELeighteen$ is bounded by its initial value. The non-linear estimates \eqref{enNLPersistanceOfRegularity}-\eqref{enNLMorawetz} are not needed in the linear case. We note the loss of angular regularity comes entirely from \eqref{n2.3} as a consequence of the probable lack of angular derivative decay near the photon sphere. There is a lose of $2\varepsilon$ in the commutator estimate, and, in distributing, this between the energy and the conformal charge, we put $18$ factors on $\ELeighteen$ so that we can get none on $\ConfChg$ and get a factor of $N^{-1/4}$. 

\newstepsequence
\newstep{Persistence of regularity, \eqref{n2.1}.} 
\label{nStep1}
For this step, $\alpha$ will be an index of angular regularity, and not the value of $r$ corresponding to the peak of the effective potential or to the photon sphere. We apply the operator $L^\alpha$ to the wave equation,
\begin{align}
\dtt(L^\alpha\phi)
=\drr(L^\alpha\phi) - V(L^\alpha\phi) -V_L(-\wlap)(L^\alpha\phi) -L^\alpha(fF'\phi) ,
\hidenumber \end{align}
and consider the densities and relations with general nonlinearity $\mathcal{N}=L^\alpha(fF'\phi)$, 
\begin{align}
e_0[v,w] =& \frac12( w^2 +(\dr v)^2 + Vv^2 +V_L|\nabla_{\cpctW}v|^2) ,&
p_{0*}[v,w] =& w\dr v ,&
\vec{p}_{0\omega} =& w \nabla_{\cpctW}v .
\hidenumber \end{align}
If $\psi$ solves the general nonlinear wave equation, 
\begin{align}
\dtt(\psi)
=\drr(\psi) - V(\psi) -V_L(-\Delta_{S^2})(\psi) -\mathcal{N} ,
\hidenumber \end{align}
then there is a relation for $e_0$, and a formula for the time derivative of its spatial integral. 
\begin{align}
\dt e_0[\psi,\dot\psi]
=& \dr p_{0*}[\psi,\dot\psi] + \nabla_{\cpctW}\cdot\vec{p}_{0\omega}[\psi,\dot\psi] - \dot\psi\mathcal{N} ,&
\Dt E_0[\psi,\dot\psi](t)
=& \Dt \int_\starman e_0[\psi,\dot\psi] \dThree 
=- \int_\starman \dot\psi\mathcal{N} .
\hidenumber \end{align}
Taking $\psi=L^\alpha\phi$, 
\begin{align}
\ELalpha(t)
=& E_0[\psi,\dot\psi](t)
=\int_0^t\int_\starman L^\alpha\dot\phi L^\alpha fF'(|\phi|^2)\phi \dThree d\tau +\ELalpha(0) \hidenumber\\
\leq& \int_0^t \|L^\alpha \dot\phi\| \|L^\alpha fF'(|\phi|^2)\phi\| d\tau +\ELalpha(0) .
\hidenumber \end{align}
We now take $\alpha=18\varepsilon$. 

\newstep{Temporal refinement, \eqref{n2.3}.} 
\label{nStep2}
We consider the operator $L^{2\varepsilon}\Gammacrown$ and apply the Heisenberg identity
\begin{align}
\Dt(\langle \dot\phi,L^{2\varepsilon}\Gammacrown\phi\rangle - \langle\phi,L^{2\varepsilon}\Gammacrown\dot\phi\rangle)
=& \langle \phi,[\linH,L^{2\varepsilon}\Gammacrown]\phi\rangle+\langle \dot\phi,L^{2\varepsilon}\dot\Gammacrown\phi\rangle - \langle\phi,L^{2\varepsilon}\dot\Gammacrown\dot\phi\rangle \hidenumber\\
&+\langle fF'\phi,L^{2\varepsilon}\Gammacrown\phi\rangle - \langle\phi,L^{2\varepsilon}\Gammacrown fF'\phi\rangle
\label{n5.1}
\end{align}
Ignoring factors of $L^\varepsilon$, our method is to show that the term on the left, under the time derivative is conformally bounded and that the linear term on the right dominate $tL^2\chia$ in expectation value, so that the time integral of the expectation value of $tL^2$ is conformally bounded. We now show this with the additional factors of $L^\varepsilon$ as required. The last two terms in equation \eqref{n5.1} are not present in the linear case. 

We begin by bounding the skew expectation at the end points. Since $\Gammacrown=(1+t)(\Gammahat+\CGamma\chit\gamma\chit)$, and $\Gammacrown=\Gammacrown\chit$, we can estimate the norm of $\Gammacrown$, $\|\Gammacrown\phi\|^2=\|\Gammacrown\chit\phi\|^2\leq (1+t)E[\chit\phi,\chit\dot\phi]\leq\ConfChg[\phi,\dot\phi]$. Since $L^{2\varepsilon}\Gammacrown$ is anti-self-adjoint and additional light cone localization can be applied to $\Gammacrown$, 
\begin{align}
|\langle\dot\phi,L^{2\varepsilon}\Gammacrown\phi\rangle-\langle\dot\phi,L^{2\varepsilon}\Gammacrown\phi\rangle|
=&2|\langle L^{2\varepsilon}\chit\dot\phi,\Gammacrown\phi\rangle| \hidenumber\\
\leq& C \ELtwo^\frac12 \ConfChgphi^\frac12 .
\label{n6.1}
\end{align}

We proceed temporarily using the notation $\OSterbenzEpsilon$ to denote $\OSterbenz$ with respect to $L^\varepsilon\phi$ instead of $\phi$. From lemma \ref{lGammacrownMajoratesLOneMinusEpsilon} applied to $L^\varepsilon\phi$, since $L^\varepsilon$ commutes with $\Gammacrown$ and $\linH$, 
\begin{align}
\langle \phi,[\linH,L^{2\varepsilon}\Gammacrown]\phi\rangle+\langle \dot\phi,L^{2\varepsilon}\dot\Gammacrown\phi\rangle - \langle\phi,L^{2\varepsilon}\dot\Gammacrown\dot\phi\rangle 
\geq& C\langle \phi, tL^2\chia\phi\rangle + \OSterbenzEpsilon
\label{n6.2}
\end{align}

The nonlinear terms are 
\begin{align}
|\langle fF'\phi,L^{2\varepsilon}\Gammacrown\phi\rangle - \langle\phi,L^{2\varepsilon}\Gammacrown fF'\phi\rangle|
\leq& 2|\langle\Gammacrown\phi,\chit L^{2\varepsilon} f F'(|\phi|^2)\phi \rangle| .
\label{neNLPartOfMorawetz}
\end{align}

We integrate the Heisenberg relation \eqref{n5.1} and substitute \eqref{n6.1}-\eqref{neNLPartOfMorawetz}, 
\begin{align}
C \ELtwo^\frac12 \ConfChgphi^\frac12
\geq& \int_0^t \langle\phi,\tau L^2\chia\phi\rangle + \OSterbenzEpsilon d\tau\hidenumber\\
& -C |\int_0^t \langle fF'\phi,L^{2\varepsilon}\Gammacrown \phi\rangle| d\tau
\hidenumber \end{align}
Solving for the time integral of $tL^2$, we get
\begin{align}
\int_0^t \langle\phi,\tau L^2\chia\phi\rangle d\tau
\leq& C \supt\ELtwo(\tau)^\frac12 \supt\ConfChgphi(\tau)^\frac12 \hidenumber\\
&+ \int_0^t \OSterbenzEpsilon d\tau \hidenumber\\
&+C |\int_0^t \langle fF'\phi,L^{2\varepsilon}\Gammacrown \phi\rangle d\tau|. 
\label{n7.2}
\end{align}
To prove \eqref{n2.3}, it remains to bound the integral on the right. We decompose the integral onto harmonics and then decompose into small and large time relative to the harmonic parameter. This allows us to make estimates based on the supremum of norms of $\phi$ rather than the sum of suprema of the projections of $\phi$. We recall that $\Ptl$ denotes the projection onto the $l^\text{th}$ harmonic. We use $\lambda^2$ to refer to the eigenvalue of $1-\wlap$. 
\begin{align}
\int_0^t \OSterbenzEpsilon d\tau
\leq& C\sum_l \lambda^{2\varepsilon} \int_0^{N\lambda^{8\varepsilon}} E[\Ptl\phi,\Ptl \dot\phi] d\tau \hidenumber\\
&+ C\sum_\lambda \lambda^{2\varepsilon} \int_{N\lambda^{8\varepsilon}}^t \ConfChg[\Ptl\phi,\Ptl \dot\phi] \tau^{-2} d\tau . \hidenumber\\
\leq& C\sum_\lambda N^2 \lambda^{18\varepsilon} \int_0^{N\lambda^{8\varepsilon}} E[\Ptl\phi,\Ptl \dot\phi](\tau) (1+\tau)^{-2} d\tau 
+ C\sum_\lambda N^{-\frac14} \int_{N\lambda^{8\varepsilon}}^t \ConfChg[\Ptl\phi,\Ptl \dot\phi] \tau^{-\frac74} d\tau \hidenumber\\
\leq& C N^2 \int_0^{t} \ELeighteen(\tau) (1+\tau)^{-2} d\tau 
+ C N^{-\frac14} \int_{0}^t \ConfChg[\phi,\dot\phi] \tau^{-\frac74} d\tau \hidenumber\\
\leq& C N^2 \supt \ELeighteen(\tau) \int_0^t (1+\tau)^{-2} d\tau 
+ N^{-\frac14} \supt \ConfChgphi(\tau) \int_0^t (1+\tau)^{-\frac74} d\tau \hidenumber\\
\leq& C N^2 \supt \ELeighteen(\tau) 
+ CN^{-\frac14} \supt \ConfChgphi(\tau) 
\label{n8.1}
\end{align}

The contribution for $t<N\lambda^{8\varepsilon}$ can be reduced to $N \ELten$ in the linear case, but there's no major gain in the argument. From \eqref{n7.2} and \eqref{n8.1}, 
\begin{align}
\int_0^t \langle\phi,L^2\chia\phi\rangle d\tau
\leq& C \supt \ELtwo(\tau)^\frac12 \supt \ConfChgphi(\tau)^\frac12 \hidenumber\\
&+ C N^2 \supt \ELeighteen(\tau) + CN^{-\frac14} \supt\ConfChgphi(\tau)\hidenumber\\
& +C |\int_0^t \langle fF'\phi,L^{2\varepsilon}\Gammacrown \phi\rangle d\tau| \hidenumber\\
\label{n9.1} 
\leq& C N^2 \supt \ELeighteen(\tau) + CN^{-\frac14} \supt\ConfChgphi(\tau)
+C |\int_1^t \langle fF'\phi,L^{2\varepsilon}\Gammacrown \phi\rangle d\tau| 
\end{align}

\newstep{NL persistence of regularity, \eqref{enNLPersistanceOfRegularity}.}
We continue working with an arbitrary angular regularity exponent $\alpha$. We begin by stating a fractional chain rule estimate in $\Reals^{n}$ \cite[Proposition 5.1, p112]{Taylor}. In $\Reals^n$, with $s\in(0,1)$, $1/q_0=1/q_1+1/q_2$, and $x\mapsto|x|^{p-1}x$ a $C^1$ function, 
\begin{align}
\| |u|^{p-1}u \|_{H^{s,q_0}} \leq C \| |u|^{p-1} \|_{\Lebesgue^{q_1}} \| u\|_{H^{s,q_2}} .
\end{align}
Here $H^{s,q}$ denotes the space of functions with $s$ derivatives in $L^q$, where $s$ derivatives are defined by the Fourier multiplier $(1+|\xi|^2)^{s/2}$. By taking a partition of unity, we have the same result on the sphere. Since $L$ corresponds to the Bessel potential $(1+|\xi|^2)^{1/2}=(1-\Delta)^{1/2}$, taking $q_0=2$, $q_1=3$, and $q_2=6$, and then integrating in space against the weight $f$ we have 
\begin{align}
\| L^{18\varepsilon} f(r_*) F'(|\phi|^2)\phi\|
\leq& \| f V_L^{-\frac13} |\phi|^{p-1} \|_3 \| V_L^\frac13 L^{18\varepsilon} \phi \|_6 . \hidenumber 
\end{align}

Under condition \ref{CNLrBetter}, the $\Lebesgue^6$ norm of $L^{18\varepsilon}\phi$ can be estimated by the Sobolev estimate \eqref{eSobolev} with the non-derivative terms estimated by the Hardy inequality \eqref{eHardyPoincareSterbenz}. For $2<3(p-1)<6$, the second can be estimated by the the conformal-Sobolev estimates \eqref{enpKlainermanSobolev}. 
\begin{align}
\| L^{18\varepsilon} f(r_*) F'(|\phi|^2)\phi\|
\leq& \ELeighteen^\frac12 \ConfChgphi^{p-1} (1+t)^{\frac{-3(p-1)+2}{3}} \hidenumber\\
\leq& \ELeighteen \ConfChgphi^{p-1} (1+t)^{-p+\frac53} \hidenumber
\end{align}

This can now be integrated in time if $p>\frac83$, 
\begin{align}
\int_0^t \| L^{18\varepsilon} f(r_*) F'(\phi^2)\phi\| d\tau 
\leq& \int_0^t \ELeighteen \ConfChgphi^{\frac{p-1}{2}} (1+\tau)^{-p+\frac53} d\tau \hidenumber\\
\leq& \supt \ELeighteen(\tau) \supt \ConfChgphi(\tau)^{\frac{p-1}{2}} \hidenumber 
\end{align} 
It is at this stage that we have used the hypothesis $8/3<p$. 

\newstep{NL Conformal, \eqref{enNLConformal}.}
Since $F(|\phi|^2)$ is positive, and $2f+r_*f'$ is bounded by a compactly supported function, and hence by $f$ itself, 
\begin{align}
\int_\starman (2f+r_*f')F(|\phi|^2) \dThree \leq C \int_\starman f F(|\phi|^2) \dThree \leq C E. 
\hidenumber \end{align}

For this argument, we use $\tilde\chi$ to denote compactly supported, positive functions, which may vary from line to line, with $C$ or $\bndd$ for constants or bounded operators. Let $\tilde\chi$ be a function which dominates $2f+r_*f'$, 
\begin{align}
\int_\starman (2f+r_*f')F(|\phi|^2) \dThree 
\leq& \int_\starman \tilde\chi F(|\phi|^2) \dThree 
\hidenumber \end{align}
Since $F$ is a power semi-linearity, up to constants $F(|\phi|^2)=C F'(|\phi|^2)|\phi|^2$, so 
\begin{align}
\int_\starman \tilde\chi F(|\phi|^2) \dThree 
\geq& C \|\tilde\chi \phi\| \| \tilde\chi F'(|\phi|^2)\phi \| 
\hidenumber \end{align}
These can be estimated by \eqref{eL2bound} and \eqref{enpKlainermanSobolev}, 
\begin{align}
\|\tilde\chi \phi\| \|\tilde\chi F'(\phi^2)\phi \| 
\leq& C \ConfChg[\tilde\chi\phi,\tilde\chi\dot\phi]^{\frac{p+1}{2}} t^{-p-1}. \hidenumber 
\end{align}
Note that more powers of $t$ can be controlled than usual because $\tilde\chi\phi$ is compactly supported, so that $\|\tilde\chi\phi\|$ and $\|\dr(\tilde\chi\phi)\|$ both decay like $t^{-1}\ConfChg[\phi,\dot\phi]$, instead of being merely bounded like $\|\phi\|$ and $\|\dr\phi\|$. 

Integrating in time, we have
\begin{align}
\int_0^t \tau \int_\starman (2f+r_*f')F(|\phi|^2) \dThree d\tau
\leq& C \int_1^{N^{\frac14}} E \tau d\tau 
+\int_{N^{\frac14}}^t  \ConfChg^{\frac{p+1}{2}} \tau^{-p} d\tau \hidenumber\\
\leq& C N^{\frac12} E + C N^{\frac{-p+1}{4}} \supt\ConfChg(\tau)^{\frac{p+1}{2}} 
\hidenumber \end{align}

Note that, if we had merely required that the nonlinear trapping terms were merely bounded by $f$, then we would only control $\|\phi\|$ and $\|\dr\phi\|$ by $\ConfChg[\phi,\dot\phi]$. In this case, we would lose two factors of $\tau$ in the estimate, and have $\tau^{2-p}$ in the integrand. Since the Sobolev estimate requires $p\leq3$, this would not be integrable. 

\newstep{NL Morawetz, \eqref{enNLMorawetz}.}
The nonlinear terms can be estimated by using the Sobolev estimate \eqref{eSobolev} and conformal-Sobolev estimate \eqref{enpKlainermanSobolev} as in the persistence of regularity argument. 
\begin{align}
|\langle fF'(|\phi|^2)\phi,L^{2\varepsilon}\Gammacrown\phi\rangle - \langle\phi,L^{2\varepsilon}\Gammacrown fF'(|\phi|^2)\phi\rangle|
\leq& 2|\langle\Gammacrown\phi,\chit L^{2\varepsilon} f|\phi|^{p-1}\phi\rangle|\hidenumber\\
\leq& C \|\Gammacrown\phi\| \| L^{2\varepsilon}f|\phi|^{p-1}\phi\| 
\hidenumber \end{align}
We now estimate $\|\Gammacrown \phi\|$ by 
\begin{align}
\|\Gammacrown\phi\| \leq C\ConfChg[\chit\phi,\chit\dot\phi]\leq C\ConfChg[\phi,\dot\phi] , 
\hidenumber \end{align}
and $\| L^{2\varepsilon}f|\phi|^{p-1}\phi\|$ the same way as in the non linear persistence of regularity argument, and integrate all of this in time. 
\begin{align}
\| L^\alpha f(r_*) F'(\phi^2)\phi\|
\leq& \ELtwo^\frac12 \ConfChgphi^{\frac{p-1}{2}} (1+t)^{-p+\frac53} , \hidenumber\\
\int_0^t |\langle fF'(|\phi|^2)\phi,L^{2\varepsilon}\Gammacrown\phi\rangle - \langle\phi,L^{2\varepsilon}\Gammacrown fF'(|\phi|^2)\phi\rangle| d\tau
\leq& \supt \ELeighteen(\tau)^\frac12 \supt \ConfChgphi(\tau)^{\frac{p}{2}}
\hidenumber \end{align}

\newstep{Closing the boot strap.}
We summarize the results of estimates \eqref{n2.1} - \eqref{enNLConformal} as
\begin{align}
\ELeighteen(t)\leq&\ELeighteen(0) + C\supt\ELeighteen(\tau)\supt\ConfChgphi(\tau)^{\frac{p-1}{2}}
\label{n10.1} ,\\
\ConfChgphi(t)
\leq&\ConfChgphi(0) \hidenumber\\
&+ C N^2\supt\ELeighteen(\tau) + CN^{-\frac14}\supt\ConfChgphi(\tau) \hidenumber\\
&+ C\supt\ELeighteen(\tau)^\frac12 \supt\ConfChgphi^\frac{p}{2} \hidenumber\\
&+ C N^\frac12 E[\phi,\dot\phi] + CN^{-\frac{p}{4}}\supt\ConfChg{\phi}(\tau)^\frac{p+1}{2}
\hidenumber \end{align}

In terms of the norm
\begin{align}
\BSSnorm{\phi(t)}^2=&
N^\frac12 E + \supt \ConfChgphi(\tau) + N^\frac{9}{4} \supt \ELeighteen(\tau) , 
\hidenumber \end{align}
these two estimates can be combined as
\begin{align}
\BSSnorm{\phi(t)}^2
\leq& \BSSnorm{\phi(0)}^2 
+ C N^{-\frac14} \BSSnorm{\phi(t)}^2 
+ C N^{\frac14} \BSSnorm{\phi(t)}^{p+1} .\hidenumber 
\end{align}

We now take $N$ sufficiently large, so that $C N^{-1/4} < \frac12$. 
\begin{align}
\BSSnorm{\phi(t)}^2\leq& C\BSSnorm{\phi(0)}^2 + CN^{-\frac{p}{4}} \BSSnorm{\phi(t)}^{p+1} . \hidenumber 
\end{align}
For sufficiently small initial data in the $\BSSnorm{\phi}$ norm, which means in both the $\ELeighteen$ and $\ConfChgphi$ norms, the $\BSSnorm{\phi}$ norm is uniformly bounded. As a consequence, the $\ELeighteen$ and $\ConfChgphi$ norms are also uniformly bounded. 
\end{proof}

From the boundedness of the conformal charge, the Sobolev estimate \eqref{eSobolev} and the control on the $L^2$ norm by conformal charge in equation \eqref{eL2bound} in the proof of \eqref{lpKlainermanSobolev}, the $L^6$ norm decays
\begin{align}
\int_\starman V_L^2 |\phi|^6 \dThree \leq \ConfChg^3 (1+t^2)^{-2} .
\hidenumber \end{align}
Because of the time reversibility of the wave equation, the same estimate holds as $t\rightarrow-\infty$. 

For the wave equation on the warp product manifolds considered, this means that
\begin{align}
\int_\starman |\tilde\phi|^4 \dTildeThree
=&\int_\starman |\phi|^4 r^{-2} \dThree \hidenumber\\
\leq&(\int_\starman |\phi|^2 \dThree)^\frac12 (\int_\starman V_L^2 |\phi|^6 \dThree)^\frac12 \hidenumber\\
\leq& \ConfChg (1+t^{2})^{-1} . \hidenumber 
\end{align}
Since $t^{-2}$ is time integrable and $\ConfChg$ is bounded, the spatial integral of $|\tilde\phi|^4$ is time integrable. Thus, $|\tilde\phi|^4$ is integrable in space-time. 

For the wave equation on the Schwarzschild solution, the situation is a little more complicated. We consider the space-time $L^4$ norm directly, and recall that the natural measure has an additional factor of $(1-2M/r)$. 
\begin{align} 
\int_{-\infty}^\infty \int_\starman |\tilde\phi|^4 r^2 dr d^2\omega dt
=& \int_{-\infty}^{\infty} \int_\starman (1-2M/r)r^{-2} |\phi|^4 \dThree dt \hidenumber\\
\leq& \int_{-\infty}^{\infty} (\int_\starman |\phi|^2 \dThree)^\frac12 (\int_\starman V_L^2 |\phi|^6 \dThree)^\frac12 dt\hidenumber\\
\leq& \int_{-\infty}^{\infty} \ConfChg[\phi(t),\dot\phi(t)]^2 (1+t^2)^{-1} dt \hidenumber\\
& \leq C (\ConfChg[\phi(0),\dot\phi(0)] + E[L^\varepsilon\phi(0),L^\varepsilon\dot\phi(0)])^2 . \hidenumber 
\end{align} 

\subsection{Summary of estimates}
We now summarize the key estimates for linear solutions. These hold both for wave equation on the Schwarzschild solution and on Riemannian manifolds under the conditions in section \ref{sSetup}. 
\begin{align}
\int \|(1+r_*^2)^{-1} \phi\|^2 d\tau \leq& C E \\
\int_\starman V_L |\nabla_{\cpctW} \phi|^2 \dThree d\tau \leq& C \ConfChg(t) (1+t^2)^{-1} \\
\int_\starman \chi(r_*) (|\dot\phi|^2 + |\phi'|^2 +|\phi|^2 + |\nabla_{\cpctW} \phi|^2) \dThree \leq& C \ConfChg(t) (1+t^2)^{-1} \\
\int_\starman V_L^2 |\phi|^6 \dThree \leq& C \ConfChg(t)^3 (1+t^2)^{-2} 
\end{align}
where $\chi(r_*)$ refers to any compactly supported function and
\begin{align}
E =& E[\phi(0),\dot\phi(0)] \\
\ConfChg(t)\leq&C (\ConfChg[\phi(0),\dot\phi(0)] + E[L^\varepsilon\phi(0),L^\varepsilon\dot\phi(0)]) 
\end{align}
For small-data, non-linear problems, there is a correction term in the first estimate and the others remain unchanged. 

From these, solutions to the original, geometrically defined wave equation $\dalembertian\tilde\phi=0$ satisfy the following space-time integral estimate, in terms of the geometrically defined $4$-volume. 
\begin{align}
\int |\tilde{\phi}|^4 d^4\text{vol} \leq& C (\ConfChg[\phi(0),\dot\phi(0)] + E[L^\varepsilon\phi(0),L^\varepsilon\dot\phi(0)])^2
\end{align}
The region of integration is the exterior region, region $\text{I}$, on the Schwarzschild solution or the entire space-time $\Reals\times\starman$ on a Riemannian manifold. On a Riemmanian manifold, for the semi-linear equation $\dalembertian\tilde\phi=|\phi|^{p-1}\phi$ with small initial data and $p\in(8/3,3)$, the same estimate holds. 

\appendix
\section{Estimates on the partial phase space observables}
\label{saPhaseSpaceEst}

We present three results from our previous work \cite{BlueSoffer3}. These are phase space localized estimate which we use to build the estimates involving the phase space observable $\Gammahat$ and reduce the loss of angular regularity at the photon sphere to $\epsilon$ derivatives. The first two estimates state that the partial phase space observables majorate powers of $L$, and the third, that the partial phase space observables are energy bounded. For these recall that $\Gammanm$ is defined in terms of $\delta$ and $\epsilon$. 

\begin{lemma}
\label{lEstimatingThePhaseSpaceCommutatorBeneathStrip}
For \psianicefn, $\delta, \epsilon >0$, and $n\in[0,\frac12]$, there is a constant $C$ such that 
\begin{align}
\langle \psi,[\linH,\Gammanhalf]\psi\rangle\geq& C\|L^{\frac{3/2+n-\epsilon}{2}}\pnn{\pn}\chia \psi\|^2-\ErrorTermsnPsi
\hidenumber \end{align}
\end{lemma}

\begin{lemma}
\label{lEstimatingThePhaseSpaceCommutatorInStrip}
For \psianicefn, $\delta, \epsilon >0$, and $n\in[0,\frac12]$,
\begin{align}
\langle \psi,[\linH,\Gammann]\psi\rangle\geq& C\|L^{1-\frac\epsilon2-\delta}\pni{\pn}\chia \psi\|^2 -\ErrorTermsnPsi
\hidenumber \end{align}
\end{lemma}

\begin{lemma}
\label{lGammanmIsEBndd}
For \psianicefn{} and $n\leq m\leq1/2$, 
$\Gammanm$ is energy bounded: 
\begin{align}
C_{nm}\|\Gammanm\psi\|^2 \leq E[\psi,\dot\psi] .
\hidenumber \end{align}
\end{lemma}

In section \ref{sAngularRefinements}, we define the operators $L$, $\qm$, and $\pn$; the concept of lower order, $\ErrorTermsnPsi$; the phase space localizations, $\pnna{\pn}$, $\pnn{\pn}$, $\pni{\pn}$; and the partial phase space observables, $\Gammanm$.

Lemmas \ref{lEstimatingThePhaseSpaceCommutatorBeneathStrip}-\ref{lGammanmIsEBndd} are essentially phase space localized versions of the estimates on $\gamma$ from section \ref{sLocalDecay}. Since the partial phase space observables are phase space localized versions of $\gamma$, with extra powers of $L$, the results would follow immediately if phase space localization could be commuted through spatial localization freely. The essence of the proofs in this section is that commuting the phase space localization terms only generates lower order terms. We state the general commutator expansion theorem of Sigal-Soffer \cite{Soffer} and specialize to the case $[F_1(\qm),F_2(\pn)]=L^{m+n-1}\bndd$. 

We use the notation $\Fourier{\cdot}$ to denote the Fourier transform. The $k$th commutator, $\ad_A^k(B)$, is the quadratic form $\langle \ad_A^{k-1}(B)\psi,A\psi\rangle-\langle A\psi,\ad_A^{k-1}(B)\psi\rangle$. If $\ad_A^{k}$ extends to a bounded operator, the same notation refers to this operator. 

\begin{theorem}[Commutator Expansion Theorem]
\label{tIIp2.1}
If $n>0$ is an integer, $A$ is a self-adjoint operator, $F_1$ is a self-adjoint operator for which for $1\leq k\leq n$, $\ad_A^k(F_1)$ extends to a bounded operator, and $F_2(x)$ is a smooth function satisfying $\|\Fourier{F_2^{[n]}}\|_1\leq\infty$, where $\Fourier{\cdot}$ denotes the Fourier transform, 
then if $[F_1,F_2(A)]$ is defined as a form on the domain of $A^n$, 
\begin{align}
\label{IIp2.2}
[F_1,F_2(A)]=\sum_{k=1}^{n-1}\frac{1}{k!} F_2^{[k]}(A) \ad_A^k(F_1) +R_n
\end{align}
in the form sense with the remainder $R_n$ satisfying
\begin{align}
\|R_n\|_\opnorm\leq& C\|\Fourier{F_2^{[n]}}\|_1 \|\ad_A^n(F_1)\|_\opnorm
\hidenumber \end{align}
Consequently, $[F_1,F_2(A)]$ defines an operator on the domain of $A^{n-1}$. 
\end{theorem}

We now specialize to the case that $F_1$ is a function of $\qm$ and $A=\pn$. The commutator $[\qm,\pn]$ can be computed explicitly as $iL^{m+n-1}$ which has a negative power of $L$ for $n\leq m<\frac12$. 

\begin{lemma}
Let $k$ be a positive integer, $F_1(x)$ be a smooth function with, for $1\leq j\leq k$, $\|F_1^{[j]}\|_\infty<\infty$, and $F_2(x)$ be a smooth function with $\|\Fourier{F_2^{[k]}}\|_1\leq\infty$, then there is a constant $C_k$ depending only on $k$, such that
\begin{align}
L^{k(1-m-n)}[F_1(\qm),F_2(\pn)]
=& \sum_{j=1}^{k-1} i^j L^{(k-j)(1-m-n)} F_1^{[j]}(\qm)F_2^{[j]}(\pn) + R_k \hidenumber\\
\|R_k\|\leq& C_k \| F_1^{[k]}\|_\infty \|\Fourier{F_2^{[k]}}\|_1 .
\hidenumber \end{align}
\end{lemma}
\begin{proof}
This follows from making a harmonic decomposition, where $L=\lambda$, and noting the adjoint is 
\begin{align}
\ad_{\pn}^j(F_1(\qm))
=&i\lambda^{n-1}\dr \ad_{\pn}^{j-1}(F_1(\qm)) \hidenumber\\
=&i^j\lambda^{j(n-1)} (\dr)^jF_1(\lambda^mr_*) \hidenumber\\
=&i^j \lambda^{j(m+n-1)} F_1^{[j]}(\lambda^mr_*) ,
\hidenumber \end{align}
and applying the commutator expansion theorem. 
\end{proof}

The operators $\gm$, $\pnna{\pn}$ and $\pnn{\pn}$ and defined in section \ref{sAngularRefinements}. We also introduce the following functions.
\begin{align}
\qmn{x} =& \frac{b}{(1+b|x|)^2} ,&
\qmnt{x} =& \frac{1}{1+x^2} ,&
\qmf{x} =& x \sqrt{\frac{g(x)}{x}} ,\hidenumber\\
\pnnb{\xi} =& \xi \pnnaWITHADERIV{\xi} ,&
\pnnc{\xi} =& \sqrt{\pnna{\xi}(\pnna{\xi}+2\xi\pnnaWITHADERIV{\xi})} .
\hidenumber \end{align}
The operators $\qmn{\qm}$, $\qmnt{\qm}$, $\qmf{\qm}$, $\pnnb{\pn}$, and $\pnnc{\pn}$ are defined by the spectral theorem. Two important properties of these functions are 
\begin{align}
\label{ps3.2}
\qmnt{x}^2+\qmf{x}^2 >&C ,\\
\label{ps3.2b}
\pnnc{\xi}>&C\pnna{\xi}
\end{align}

In applying the order reduction lemma, it is useful to note that, if $F_1$ is a $C^1$ function with bounded derivative, and if $\pnnfn(\xi)$ is one of $xi$, $\pnna{\xi}$, $\pnnb{\xi}$, or $\pnnc{\xi}$, then
\begin{align}
\label{ps3.3}
[F_1(x),\pnnfn(\xi)]=& L^{n+m-1}\bndd_1 \hidenumber\\
\bndd_1 \leq& C \|F_1'\|_\infty ,
\end{align}
where the notation $\bndd$ is used to denote an arbitrary bounded operator which commutes with $L$, the same way that $C$ is use to denote an arbitrary constant. In addition, if $F_1$ has a bounded second derivative, then for $\pnnfn(\xi)=\xi\pnna{\xi}$, 
\begin{align}
\label{ps3.4}
[F_1(\qm),\pnnfn(\pn)]=& L^{n+m-1}\bndd \hidenumber\\
\| \bndd\| \leq C (\|F_1'\|_\infty + \|F_1''\|_\infty) .
\end{align}

In particular, equations \eqref{ps3.3} and \eqref{ps3.4} hold for $F_1$ being any Schwartz class function, $\qmnt{x}$, or $\qmf{x}$. These estimates follow from the order reduction lemma and the fact that the Fourier transform of $\pnnfn'$ is in $L^1$ because the high degree of regularity for $\pnnfn$ translates into decay properties for the Fourier transform. For $\pnnfn(x)=x\pnna{x}$, $\pnnfn''$ has Fourier transform in $L^1$ and $\pnnfn'$ in $L^\infty$. 

We now prove lemma \ref{lGammanmIsEBndd}. The method of commuting localization to get lower order terms will be typical of later arguments. We will also use the notation $Q_1\lesssim Q_2$ to denote there is a constant $C$ such that $Q_1\leq C Q_2$. 

\begin{proof}[Proof of lemma \ref{lGammanmIsEBndd}]
We start by dropping a spatial localization and expanding  $\gammam$
\begin{align}
\| \Gammanm \psi\|^2
=& \| \chia \pnnan L^{n-\epsilon} \gammam \pnnan \chia \psi\|^2 \hidenumber\\
\lesssim& \| \pnnan \dr \gm L^{n-\epsilon} \pnnan \chia \psi \|^2
+ \|L^m \qmn{\qm} L^{n-\epsilon} \pnnan \chia \psi\|^2 .
\hidenumber \end{align}
We rewrite the first term in terms of $\pn$ and drop a bounded operator in the second, leaving the bound
\begin{align}
\| L^{1-\epsilon}\pnnan \pn \gm \pnnan \chia \psi \|^2 + \| L^{m+n-\epsilon} 
\chia \psi \|^2 . 
\hidenumber \end{align}
We commute phase space localizations $\pnnan \pn$ and $\gm$ in the first term and write the commutator at the end of the expression, as we will typically do. The bound becomes
\begin{align}
\| L^{1-\epsilon} \gm \pnnan^2\pn \chia\chi\psi\|^2 + \| L^{m+n-\epsilon} \chia \psi\|^2
+ \|L^{1-\epsilon} [\pnnan\pn,\gm]\pnnan\chia\psi\|^2 .
\hidenumber \end{align}
We now drop the bounded operator $\gm$ in the first term and apply the order reduction lemma in the last to get the upper bound
\begin{align}
\|L^{1-\epsilon} \pnnan^2\pn \chia\psi\|^2 + \| L^{m+n-\epsilon}\chia\psi\|^2 + \| L^{m+n-\epsilon} \bndd \pnnan \chia \psi \|^2 .
\hidenumber \end{align}
From the decay of $\pnnan$, we have that 
\begin{align}
\pnna{\xi}^2 \xi \leq |\xi|^\epsilon .
\hidenumber \end{align}
Applying $f^{1-\epsilon} g^\epsilon \leq |f| + |g|$, and interpolating $L^{m+n}$ between $L$ and $1$, gives
\begin{align}
\|\Gammanm\psi\|^2 
\lesssim& \|L\chia\psi\|^2 + \|\dr\chia\psi\|^2 + \|L\chia\psi\|^2 + \|\chia\psi\|^2 \hidenumber\\
\lesssim& E[\psi,\dot\psi] .
\hidenumber \end{align}
\end{proof}

To prove lemmas \ref{lEstimatingThePhaseSpaceCommutatorBeneathStrip} and \ref{lEstimatingThePhaseSpaceCommutatorInStrip}, we prove preliminary bounds on the commutators involving $\linH$. These bounds are localized in both $\qm$ and $\pn$. The key idea, as demonstrated in the previous proof, is that commutators are lower order terms. 

\begin{lemma}
For \psianicefn, $n<1/2$, $n\leq m\leq 1/2$, and $\epsilon, \delta > 0$, there are constants depending only on $n$, $m$, $\epsilon$, $\delta$, such that
\begin{align}
\label{ps5.5}
\langle \psi,[H_2,\Gammanm] \psi \rangle
=&\ErrorTermsnPsi \\
\label{ps5.6}
\langle \psi, [H_1+H_3,\Gammanm] \psi \rangle
\geq& C \langle \pnnan\chia \psi,L^{n-\epsilon}(-\dr \qmnt{\qm} \dr - L^2 \gm V_L') \pnnan{\pn}\chia \psi\rangle \\
&+\ErrorTermsnPsi .
\hidenumber \end{align}
\end{lemma}

\begin{proof}
We first consider the commutator with $H_2$. The commutator with $H_2$ is expanded by the Leibniz rule. The commutator with $\chia$ is zero.  
\begin{align}
[H_2,\Gammanm]
=&\chia[V,\pnnan]\gammam L^{n-\epsilon}\pnnan\chia + \chia\pnnan\gammam L^{n-\epsilon}[V,\pnnan]\chia \hidenumber\\
&+\chia\pnnan[V,\gammam]L^{n-\epsilon}\pnnan\chia .
\hidenumber \end{align}
The first two terms are, up to signs, adjoint. The commutators involving $\pnnan$ can be expressed using the order reduction lemma, and the commutator with $\gammam$ can be expanded as in section \ref{sLocalDecay}
\begin{align}
|\langle\psi,[H_2,\gammam]\psi\rangle|
\lesssim& |\langle\gamma\pnnan\chia\psi,L^{n-\epsilon}[V,\pnnan]\chia\psi\rangle| \hidenumber\\
&+|\langle\pnnan\chia\psi,[V,\gammam]L^{n-\epsilon}\pnnan\chia\psi\rangle| \hidenumber\\
\lesssim& \|\gamma\pnnan\chia\psi\|^2 + \| L^{n-\epsilon} L^{m-1}\bndd\chia\psi\|^2\hidenumber\\
&+ \|\pnnan\chia\psi\|^2 + \|\gm V' L^{n-\epsilon}\pnnan\chia\psi\|^2\hidenumber\\
\lesssim& \ErrorTermsnPsi .
\hidenumber \end{align}

By an identical argument, the difference between the commutator $[H_3,\Gammanm]=[(L^2-1)V_L,\Gammanm]$ and $[L^2V_L,\Gammanm]$ is lower order, and we will ignore this difference in the following arguments. 

The terms coming from $H_1+H_3$ can also be expanded using the Leibniz rule. 
\begin{align}
[H_1+H_3,\Gammanm]
=& [H_1+H_3,\chia]\pnnan L^{n-\epsilon}\gammam\pnnan\chia\hidenumber\\
&+ \chia\pnnan L^{n-\epsilon}\gammam\pnnan[H_1+H_3,\chia]\hidenumber\\
&+ \chia[H_1+H_3,\pnnan] L^{n-\epsilon}\gammam\pnnan\chia\hidenumber\\
&+ \chia\pnnan L^{n-\epsilon}\gammam[H_1+H_3,\pnnan]\chia\hidenumber\\
&+ \chia\pnnan L^{n-\epsilon}[H_1+H_3,\gammam]\pnnan\chia
\hidenumber \end{align}

The first two terms are estimated first. Since $H_3=(L^2-1)V_L$ commutes with functions of $r_*$, the contribution from $H_3$ can be ignored. The terms are, up to signs, each other's adjoints, so it is sufficient to consider the second. We work in expectation value and begin by expanding $\gammam$, 
\begin{align}
|\langle\psi,\chia\pnnan L^{n-\epsilon}&\gammam\pnnan[-\drr,\chia]\psi\rangle|\hidenumber\\
\leq& |\langle\psi,\chia\pnnan L^{n-\epsilon}\dr\gm\pnnan[-\drr,\chia]\psi\rangle|\hidenumber\\
&+ |\langle\psi,\chia\pnnan L^{n-\epsilon}L^m \qmn{\qm}\pnnan[-\drr,\chia]\psi\rangle| \hidenumber\\
\leq& |\langle\dr\chia\psi, L^{n-\epsilon}\pnnan\gm\pnnan[-\drr,\chia]\psi\rangle|\hidenumber\\
&+ |\langle L^{m+n-\epsilon}\pnnan\chia\psi, \qmn{\qm}\pnnan[-\drr,\chia]\psi\rangle| .
\hidenumber \end{align}
Expanding
\begin{align}
[-\drr,\chia] = -\dr\chia'-\chia'' ,
\hidenumber \end{align}
and applying the Cauchy-Schwartz inequality gives
\begin{align}
|\langle&\psi,\chia\pnnan L^{n-\epsilon}\gammam\pnnan[-\drr,\chia]\psi\rangle|\hidenumber\\
\lesssim& \|\dr\chia\psi\|^2 + \|L^{n-\epsilon}\pnnan\gm\pnnan[-\drr,\chia]\psi\|^2 + \|L^{m+n-\epsilon}\chia\psi\|^2 \hidenumber\\
&+ \|\pnnan\qmn{\qm}\pnnan[-\drr,\chia]\psi\|^2 \hidenumber\\
\lesssim& \|L^{n-\epsilon}\pnnan\gm\pnnan[-\drr,\chia]\psi\|^2 + \|[-\drr,\chia]\psi\| + \ErrorTermsnPsi \hidenumber\\
\lesssim& \|L^{n-\epsilon}\pnnan\gm\pnnan\dr\chia'\psi\|^2 \hidenumber\\
&+ \|L^{n-\epsilon}\pnnan\gm\pnnan\chia''\psi\|^2 + \|\dr\chia'\psi\|^2 + \|\chia''\psi\|^2 + \ErrorTermsnPsi. 
\hidenumber \end{align}
The second, third, and fourth terms are clearly lower order. It remains to estimate the first, which we rewrite and then estimate by the order reduction lemma
\begin{align}
|\langle\psi,\chia\pnnan L^{n-\epsilon}\gammam\pnnan[-\drr,\chia]\psi\rangle|
\leq& \|L^{1-\epsilon}\pnnan\gm\pnnan\pn\chia'\psi\|^2 \hidenumber\\
&+ \ErrorTermsnPsi
\hidenumber \end{align}
We now apply the order reduction lemma again
using $\pnna(\xi)$
to get a bound of the form
\begin{align}
|\langle\psi,\chia\pnnan L^{n-\epsilon}\gammam\pnnan[-\drr,\chia]\psi\rangle|
\lesssim& \| L^{1-\epsilon} \gm\pnnan^2\pn \chia'\psi\|^2 + \ErrorTermsnPsi \hidenumber\\
\lesssim& \|L^{1-\epsilon} \pnnan^2\pn \chia'\psi\|^2 + \ErrorTermsnPsi\hidenumber\\
\lesssim& \|L^{1-\epsilon} |\pn|^\epsilon \chia'\psi\|^2 + \ErrorTermsnPsi\hidenumber\\
\lesssim& \|L \chia'\psi\|^2 + \| |\pn| \chia'\psi\|^2 + \ErrorTermsnPsi\hidenumber\\
\lesssim& \|L \chia'\psi\|^2 + \|L^{n-1}\dr \chia'\psi\|^2 + \ErrorTermsnPsi
\hidenumber \end{align}
The second term is clearly lower order since $L^{n-1}$ is a bounded operator. The first is also lower order since the localization $\chia'$ vanishes in a neighborhood of $r_*=0$. 

We now turn to the remaining three terms. 
\begin{align}
\label{ps8.3}
[H_1+H_3,\Gammanm]
=&\chia[H_1+H_3,\pnnan]L^{n-\epsilon}\gammam\pnnan\chia + \chia\pnnan L^{n-\epsilon}\gammam[H_1+H_3,\pnnan]\chia \hidenumber\\
&+ \chia\pnnan L^{n-\epsilon}[H_1+H_3,\gammam]\pnnan\chia .
\end{align}
Since $H_1=-\drr$ commutes with $\pnnan$, the contribution from $H_1$ can be dropped from the commutators with $\pnnan$. The contribution from the $H_1$ commutator is
\begin{align}
\chia\pnnan L^{n-\epsilon}[H_1,\gammam]\pnnan\chia
=2\chia\pnnan L^{n-\epsilon} (-\dr \gm'\dr -\frac14\gm''') \pnnan\chia
\hidenumber \end{align}

The contribution from the $H_3$ commutators is significantly more complicated than the previous pieces, and the commutator with $\pnnan$ must be expanded to second order, 
\begin{align}
[H_3,\pnnan]
=& (\pnnaWITHADERIV{\pn}[V_L,\pn]+\pnnaWITHTWODERIV{\pn}[[V_L,\pn],\pn]+R_3)(-\wlap) .
\hidenumber \end{align}
Since the difference between the contributions from $(-\wlap)$ and $L^2$ is lower order, we will ignore it. The commutator with $V_L(L^2-(-\wlap))-V_L$ alone can be treated the same way as the commutator with $H_2$. Using the second order expansion with three sub-terms, and expanding $\gammam$ as the the sum of the two sub-terms, for each of the first two terms in equation \eqref{ps8.3} gives twelve terms, plus one more from the third term in \eqref{ps8.3}, for a total of thirteen. We treat the radial derivative terms in $\gammam$ as higher order, and group the thirteen terms in terms of formal order, with the third term in \eqref{ps8.3} treated as the same order as the top part of the other two. We strip out localization in $\chia$, and get 
\begin{align}
[H_3,\Gammanm]=&\chia(Q_1+Q_2+Q_3+Q_4+Q_5+Q_6)\chia \hidenumber\\
Q_1=L^{2+n-\epsilon} (&\pnnaWITHADERIV{\pn} [V_L,\pn]\dr g_{L^m} \pnna{\pn}\hidenumber\\
& + \pnna{\pn} g_{L^m} [V_L,\dr] \pnna{\pn}\hidenumber\\
& + \pnna{\pn} g_{L^m}\dr \pnnaWITHADERIV{\pn}[V_L,\pn])\hidenumber\\
Q_2=L^{2+n-\epsilon} \frac12 (&-\pnna{\pn}'[V_L,\pn]g_{L^m}'\pnna{\pn} +\pnna{\pn} g_{L^m}'\pnna{\pn}'[V_L,\pn])\hidenumber\\
Q_3=L^{2+n-\epsilon}(&\pnnaWITHTWODERIV{\pn}[[V_L,\pn],\pn]\dr g_{L^m}\pnna{\pn}\hidenumber\\
&+\pnna{\pn} g_{L^m}\dr\pnnaWITHTWODERIV{\pn}[[V_L,\pn],\pn])\hidenumber\\
Q_4=L^{2+n-\epsilon}\frac12(&-\pnnaWITHTWODERIV{\pn}[[V_L,\pn],\pn]g_{L^m}'\pnna{\pn}\hidenumber\\
& +\pnna{\pn} g_{L^m}'\pnnaWITHTWODERIV{\pn}[[V_L,\pn],\pn])\hidenumber\\
Q_5=L^{2+n-\epsilon}(&R_3 \dr g_{L^m}\pnna{\pn}+\pnna{\pn} g_{L^m}\dr R_3)\hidenumber\\
Q_6=L^{2+n-\epsilon}\frac12(&-R_3 g_{L^m}'\pnna{\pn}+\pnna{\pn} g_{L^m}' R_3)
\hidenumber \end{align} 

The top term $Q_1$ gives the desired estimate. All the others are lower order. We must expand to second order so that the estimates on the remainder, which can not utilize the structure of the localization, are lower order in equations \eqref{ps12.3} and \eqref{ps12.4}. 

We rearrange the terms in $Q_1$. It is at this point that the positivity of $\epsilon$ and the smoothness of the localizations in $\pn$ are used. 
\begin{align}
L^{-(2+n-\epsilon)}Q_1 
=& \pnnaWITHADERIV{\pn}(-V_L')(-iL^{n-1})\dr\gm\pnnan + \pnnan(-V_L'iL^{n-1}\gm)\pnnan\hidenumber\\
&+ \pnnan\gm\dr\pnnaWITHADERIV{\pn}(-V_L')
\hidenumber \end{align}
We rewrite this in terms of $\pn$ and commute to gather $-V_L'$ with $\gm$
\begin{align}
L^{-(2+n-\epsilon)}Q_1 
=& \pnnaWITHADERIV{\pn}(-V_L')\pn\gm\pnnan + \pnnan(-V_L'\gm)\pnnan + \pnnan\gm\pn \pnnaWITHADERIV{\pn}(-V_L') \hidenumber\\
=& \pnnaWITHADERIV{\pn}\pn(-V_L'\gm)\pnnan + \pnnan(-V_L'\gm)\pnnan + \pnnan(-V_L'\gm)\pnnaWITHADERIV{\pn}\pn \hidenumber\\
&+ \pnnaWITHADERIV{\pn} L^{n-1}\bndd\gm\pnnan - \pnnan\gm L^{n-1}\bndd . \label{ps9.3}
\end{align}
We will track inverse powers of $L$, but otherwise group bounded operators together. We now introduce the following decomposition to commute $-V_L'\gm$ through $\pn$ localization. 
\begin{align}
-V_L'\gm=& (-r_*^{-1}V_L') L^{-m} L^mr_*\gm\hidenumber\\
=&(-r_*^{-1}V_L')L^{-m}\qmf{\qm}^2
\hidenumber \end{align}
By Taylor's theorem, $-r_*^{-1}V_L'$ is an analytic function since $V_L$ is, and $V_L'$ vanishes at $r_*=0$. For a localization $\pnnfn$ being one of $\pnnaNOARG$, $\pnnbNOARG$, or $\pnncNOARG$, we write
\begin{align}
[\pnnfn(\pn),-V_L'\gm]
=&[\pnnfn(\pn),L^{-m}\qmf{\qn}^2 (-r_*^{-1}V_L')]\hidenumber\\
=&[\pnnfn(\pn),\qmf{\qm}^2]L^{-m}(-r_*^{-1}V_L') + r_*\gm[\pnnfn(\pn),-r_*^{-1}V_L'\gm] \hidenumber\\
=&L^{m+n-1}\bndd L^{-m}(-r_*^{-1}V_L') + r_*\gm L^{n-1}\bndd\hidenumber\\
=&L^{n-1}\bndd + r_*L^{n-1}\bndd .
\hidenumber \end{align}
By writing the decomposition as $-\gm V_L'$, the commutator can also be expressed as $L^{n-1}\bndd + L^{n-1}\bndd r_* $. Using this to continue with \eqref{ps9.3} gives 
\begin{align}
L^{-(2+n-\epsilon)}Q_1
=&(-V_L'\gm)\pnnan(2\pnnaWITHADERIV{\pn}\pn+\pnnan)+L^{n-1}\bndd + r_*L^{n-1}\bndd\hidenumber\\
=&(-V_L'\gm)\pnncn^2 +L^{n-1}\bndd + r_*L^{n-1}\bndd\hidenumber\\
=&\pnncn(-V_L'\gm)\pnncn +L^{n-1}\bndd + r_*L^{n-1}\bndd .
\hidenumber \end{align}
Taking the expectation value with respect to $\chia\psi$, we get
\begin{align}
\langle\chia\psi,Q_1\chia\psi\rangle
=&\|L^\frac{2+n-\epsilon}{2}(-V_L'\gm)^\frac12\pnncn\chia\psi\|^2 
+\langle\chia\psi,L^{1+2n-\epsilon}B\chia\psi\rangle+\langle r_*\chia\psi,L^{1+2n-\epsilon}B\chia\psi\rangle \hidenumber\\
=&\|L^\frac{2+n-\epsilon}{2}(-V_L'\gm)^\frac12\pnncn\chia\psi\|^2 
+ \ErrorTermsnPsi .
\hidenumber \end{align}
The localizations can be commuted, and then the $\pnncn$ localization can be replaced by $\epsilon^\frac12 \pnnan$ localization. The same arguments as before show that all the commutators are lower order. 
\begin{align}
\langle\chia\psi,Q_1\chia\psi\rangle
\gtrsim&\|L^\frac{2+n-\epsilon}{2}\pnncn(-V_L'\gm)^\frac12 \chia\psi\|^2 
+\ErrorTermsnPsi \hidenumber\\
\gtrsim&\|L^\frac{2+n-\epsilon}{2}\pnnan(-V_L'\gm)^\frac12 \chia\psi\|^2 
+\ErrorTermsnPsi .
\hidenumber \end{align}
Finally, we can commute the localizations back to get them in the standard orderer of our presentation. 
\begin{align}
\langle\chia\psi,Q_1\chia\psi\rangle
\gtrsim&\|L^\frac{2+n-\epsilon}{2}(-V_L'\gm)^\frac12\pnnan\chia\psi\|^2 
+\ErrorTermsnPsi .
\hidenumber \end{align}

We now consider the other terms in the expansion of the commutator with $H_3$. All of these are lower order. To estimate the second term, $Q_2$, we aim to absorb a factor of $L^m$ in $\qmn{\qm}$. 
\begin{align}
2L^{-(2+n-\epsilon)}Q_2
=& -\pnnaWITHADERIV{\pn}[V_L,\pn]\gmp\pnnan+\pnnan\gmp\pnnaWITHADERIV{\pn}[V_L,\pn] \hidenumber\\
=&-iL^{n-1}(-\pnnaWITHADERIV{\pn}V_L'\gmp\pnnan+\pnnan\gmp\pnnaWITHADERIV{\pn} V_L') \hidenumber\\
i2L^{-(1+2n-\epsilon)}Q_2
=&-\pnnaWITHADERIV{\pn}V_L'L^m\qmn{\qm}\pnnan+\pnnan L^m\qmn{\qm}\pnnaWITHADERIV{\pn}V_L' .
\hidenumber \end{align}
This appears to be $m$ powers of $L$ too large to be lower order. We rewrite $V_L'=r_*(r_*^{-1}V_L')$. 
\begin{align}
2L^{-(1+2n-\epsilon)}Q_2
=&-\pnnaWITHADERIV{\pn}r_*L^m(r_*^{-1}V_L')\qmn{\qm}\pnnan+\pnnan L^m\qmn{\qm}\pnnaWITHADERIV{\pn}r_*(r_*^{-1}V_L')\hidenumber\\
=&-\pnnaWITHADERIV{\pn}\qm\qmn{\qm}\pnnan(r_*^{-1}V_L')+\pnnan\qm\qmn{\qm}\pnnaWITHADERIV{\pn}(r_*^{-1}V_L') \hidenumber\\
&-\pnnaWITHADERIV{\pn}\qm\qmn{\qm}L^{n-1}\bndd + \pnnan L^m\qmn{\qm}L^{n-1}\bndd (r_*^{-1}V_L') \\
=&-\pnnaWITHADERIV\pnnan{\pn}\qm\qmn{\qm}(r_*^{-1}V_L')+\qm\qmn{\qm}\pnnan\pnnaWITHADERIV{\pn}(r_*^{-1}V_L') \hidenumber\\
&-\bndd \pnnan L^{m+n-1}L^{n-1} + \pnnaWITHADERIV{\pn}\bndd L^{m+n-1} + \bndd L^{m+n-1} . \hidenumber 
\end{align}
Since $x\qmn{x}$ is bounded and $m+n-1<0$, 
\begin{align}
\langle\chia\psi,Q_2\chia\psi\rangle = \langle\chia\psi, L^{1+2n-\epsilon}\bndd\chia\psi\rangle=\ErrorTermsnPsi . \hidenumber 
\end{align}

The third term is expanded using $[[V_L,\pn],\pn]=V_L''(-L^{2n-2})$. We aim to absorb the derivative $\dr$ into the localization function $\pnnaNOARG''$. 
\begin{align}
L^{-(2+n-\epsilon)}Q_3
=& -\pnnaWITHTWODERIV{\pn}V_L''L^{2n-2}\dr\gm\pnnan-\pnnan\gm\dr\pnnaWITHTWODERIV{\pn}V_L''L^{2n-2}\hidenumber\\
L^{-(3n-\epsilon)}Q_3
=&-\pnnaWITHTWODERIV{\pn}\dr V_L''\gm\pnnan-\pnnan\gm\pnnaWITHTWODERIV{\pn}\dr V_L'' + \bndd\hidenumber\\
L^{-(1+2n-\epsilon)}Q_3
=& -\pnnaWITHTWODERIV{\pn}\pn V_L''\gm\pnnan-\pnna{\pn}\gm\pnnaWITHTWODERIV{\pn}\pn V_L''+\bndd 
\hidenumber \end{align}
Since $\xi\pnnaWITHTWODERIV{\xi}$ is bounded, $Q_3=L^{1+2n-\epsilon}\bndd$, and $Q_3$ is lower order in expectation value relative to $\chia\psi$. 

It is clear that $Q_4$ is lower order without any rearrangement, 
\begin{align}
2Q_4
=&L^{2+n-\epsilon}(-\pnnaWITHTWODERIV{\pn}V_L''(-L^{2n-2})\gmp\pnnan+\pnnan\gmp\pnnaWITHTWODERIV{\pn}V_L''(-L^{2n-2}) )\hidenumber\\
=&L^{3n+m-\epsilon}(\pnnaWITHTWODERIV{\pn}V_L'' \qmn{\qm}\pnnan+\pnnan\qmn{\qm}\pnnaWITHTWODERIV{\pn}V_L'') .
\hidenumber \end{align}

The terms in $Q_5$ are conjugates, up to a sign, so it is sufficient to estimate one of them. 
\begin{align}
\label{ps12.3}
\langle\psi,\chia\pnnan L^{2+n-\epsilon}\dr\gm R_3\chia\psi\rangle
=&\langle\psi,\chia\pnnan L^{2+n-\epsilon}\dr\gm L^{3n-3}\bndd\chia\psi\rangle.
\end{align}
Applying Cauchy-Schwartz and dropping bounded terms,  we find this is bounded by
\begin{align}
\label{ps12.4}
\|\dr\chia\psi\|^2 + \|\gm L^{4n-1}\chia\psi\|^2 
\lesssim 
\|\dr\chia\psi\|^2 + \|L^{4n-1}\chia\psi\|^2 
=\ErrorTermsnPsi .
\end{align}

The terms in $Q_6$ are dealt with similarly
\begin{align}
\langle \psi,\chia\pnnan L^{2+n-\epsilon} \gm' R_3 \chia\psi\rangle
=& \langle\psi,\chia\pnnan L^{2+n-\epsilon} L^m \qmn{\qm} L^{3n-3}\bndd\chia\psi\rangle \hidenumber\\
=& \langle\chia\psi, L^{4n+m-1-\epsilon}\bndd\chia\psi\rangle \hidenumber\\
=&\ErrorTermsnPsi .
\hidenumber \end{align}

Gathering the result of the previous arguments, we have
\begin{align}
\langle \psi,[\linH,\Gammanm]\psi\rangle
\gtrsim&\langle \pnnan\chia \psi,L^{n-\epsilon}(-\dr\gmp\dr-\frac14\gmppp-C\gm L^2 V_L')\pnnan\chia \psi\rangle
\hidenumber \end{align}

\newcommand{\tmpAMu}{v}
This is very much like the term appearing in the Morawetz estimate, except for the scaling by $L^m$ in the function $\gm$ and the factor of $C$ on the potential term, and we proceed with a slight modification of that argument. Taking 
\begin{align}
v= L^\frac{n-\epsilon}{2} \pnnan\chia \psi
\hidenumber \end{align}
and $\chi$ to denote characteristic functions, unrelated to $\chia$, of width $O(L^{-m})$ about $r_*=0$, and multiplying by $L^{3m}$, from an analogue of the argument leading to \eqref{enForAMi}, we have
\begin{align}
\label{enInAMChiControlled}
\int L^{3m}\chi \tmpAMu^2 \dThree
=-\int L^{3m} r_*\chi' \tmpAMu^2 \dThree + 2\int L^{3m} r_*\chi \tmpAMu\dr \tmpAMu\dThree\\
\lesssim& \int L^{5m}r_*^2\chi \tmpAMu^2 \dThree
+\int L^m\chi (\dr \tmpAMu)^2 \dThree .
\end{align}
Continuing with the same argument, and multiplying the analogue of \eqref{enForAMii} by $L^{2m}$, we have, 
\begin{align}
\label{enInAMDeltaControlledByChi}
\int_{\cpctW} L^{2m}\tmpAMu(0,\omega)^2 d^2\omega
=-\int L^{2m} \chi' \tmpAMu^2 \dThree - 2\int L^{2m} \chi \tmpAMu \dr \tmpAMu \dThree \\ 
\lesssim& \int L^{2m}|\chi'| \tmpAMu^2\dThree + \int L^{3m}\chi \tmpAMu^2 \dThree + \int L^m \chi (\dr \tmpAMu)^2 \dThree
)^\frac12(\int L^{3m}r_*^2\chi \tmpAMu^2 \dThree +\int L^m\chi(\dr \tmpAMu)^2 \dThree .
\end{align}
Applying \eqref{enInAMChiControlled} to control the weighted $L^2$ norms, 
\begin{align}
\label{enInAMDeltaControlled}
L^{2m}\tmpAMu(0)^2
\leq& C_0 (\int L^{5m}r_*^2\chi \tmpAMu^2 dr_* + \int L^m\chi(\dr \tmpAMu)^2 dr_*) .
\end{align}
From integrating $\dr \gmpp \tmpAMu^2$, the analogue of \eqref{enForAMiii} is
\begin{align}
L^{3m}b^3\int \frac{\sigma(\sigma+1)}{(1+bL^m|r_*|)^{\sigma+2}} v^2 \dThree
\leq& 2\sigma b C_0\int L^{5m}br_*^2\chi \tmpAMu^2 \dThree
+(2\frac{\sigma}{\sigma+1} + 4\sigma bC_0) \int bL^m \frac{2(\dr \tmpAMu)^2}{(1+bL^m|r_*|)^\sigma} \dThree
\hidenumber \end{align}
As in \eqref{enForAMiv}, the commutator can now be estimated 
\begin{align}
\int 2\gmp(\dr \tmpAMu)^2 &- \frac12\gmppp \tmpAMu^2 -C \gm V_L' L^2 \tmpAMu^2 \dThree \hidenumber\\
\geq& \int (\frac{2}{\sigma+1} - 4\sigma bC_0) bL^m \frac{2(\dr \tmpAMu)^2}{(1+bL^m|r_*|)^\sigma} \dThree + L^{2m}b^2 \sigma \tmpAMu(0)^2 \hidenumber\\
&-\int \sigma L^{5m} b^2 C_0 r_*^2 \chi \tmpAMu^2 \dThree - C_1 \int \gm V_L' L^2 \tmpAMu^2 \dThree .
\hidenumber \end{align}
We now choose $b$ sufficiently small, independently of $L$, so that for $r_*$ sufficiently small
\begin{align}
\frac{2}{\sigma+1} - 4\sigma bC_0 \geq& \frac{1}{\sigma+1} \hidenumber\\
-C_1 \frac\epsilon3 r_*V_L' \geq&C_0 \sigma b^2 r_*^2
\hidenumber \end{align}
We now use the assumption that $m\leq1/2$. Since $\gm$ vanishes like $r_*L^m$, in the support of $\chi$, $-\frac\epsilon3 \gm L^2 V_L' > \sigma L^{5m} b^2 r_*^2\chi$, so that 
\begin{align}
\int 2\gmp(\dr \tmpAMu)^2 &- \frac12\gmppp \tmpAMu^2 -C \gm V_L' L^2 \tmpAMu^2 \dThree \hidenumber\\
\geq& \frac{b}{\sigma+1} \int L^m \frac{2(\dr \tmpAMu)^2}{(1+bL^m|r_*|)^\sigma} \dThree 
+C_1\frac{2\epsilon}{3} \int (-\gm V_L') L^2 \tmpAMu^2 \dThree .
\hidenumber \end{align}
Since $\qmn{x}$ and $\qmnt{x}$ are equivalent up to constant multiples, this gives the desired result. 
\end{proof}

The previous lemma has localization in $\qm$ and $\pn$. To prove lemmas \ref{lEstimatingThePhaseSpaceCommutatorBeneathStrip} and \ref{lEstimatingThePhaseSpaceCommutatorInStrip}, we combine the localizations $\qmn{\qm}$ and $\qmf{\qm}$, in small and large $\qm$, to get estimates localized in $\pn$. 

\begin{proof}[Proof of lemma \ref{lEstimatingThePhaseSpaceCommutatorBeneathStrip}] From the estimate in the previous lemma, \eqref{ps5.6}, we have
\begin{align}
\langle \psi,[\linH,\Gammanm]\psi\rangle
\geq& C\langle \pnnan\chia \psi,L^{n-\varepsilon}(-\dr\qmnt{\qm}\dr-\gm L^2 V_L')\pnnan\chia \psi\rangle 
\hidenumber \end{align}
We divide $\gm V_L'$ into two parts based on whether $\qm$ is smaller or larger than one. Using $\chi$ to denote a smooth characteristic function on a scale of $L^{-m}$ again, since $-\gm L^2 V_L' > L^{5m} r_*^2 \chi$, by \eqref{enInAMChiControlled}, we control $L^{3m+n-\epsilon}$ in the region $\qm < 1$. In the region $\qm>1$, $-\gm L^2 V_L' > -L^{2} \sgn(r_*) V_L'>-L^{2-m}r_*^{-1} V_L'$. We can combine both these estimates and take $m=1/2$, 
\begin{align}
\langle \psi,[\linH,\Gammanm]\psi\rangle
\geq& \langle \pnnan\chia \psi,L^{3/2+n-\epsilon}(-r_*^{-1} V_L')\pnnan\chia \psi\rangle
\hidenumber \end{align}
The factor $-r_*^{-1} V_L'$ is not really relevant, since it can be absorbed into $\chia$. To do this, we introduce a smooth, bounded function, $f$, which is equal to the inverse of $(-r_*^{-1} V_L')^{1/2}$ on the support of $\chia$, and then commute in the standard way. 
\begin{align}
\langle \psi,[\linH,\Gammanm]\psi\rangle
\geq& c\| L^{\frac{3/2+n-\epsilon}{2}} (-r_*^{-1} V_L')^\frac12 \pnnan\chia \psi\| + \ErrorTermsnPsi \hidenumber\\
\geq& c\| L^{\frac{3/2+n-\epsilon}{2}} f (-r_*^{-1} V_L')^\frac12 \pnnan\chia \psi\| + \ErrorTermsnPsi \hidenumber\\
\geq& c\| L^{\frac{3/2+n-\epsilon}{2}} \pnnan f (-r_*^{-1} V_L')^\frac12 \chia \psi\| + \ErrorTermsnPsi \hidenumber\\
\geq& c\| L^{\frac{3/2+n-\epsilon}{2}} \pnnan\chia \psi\| + \ErrorTermsnPsi 
\label{enUsefToEliminateVLp}
\end{align}
Since $\pnna{x}>c\pnn{x}$, we can replace the $\pn$ localization to get
\begin{align}
\langle \psi,[\linH,\Gammanm]\psi\rangle
\geq& c\| L^{\frac{3/2+n-\epsilon}{2}} \pnn{\pn}\chia \psi\| + \ErrorTermsnPsi 
\hidenumber \end{align}
\end{proof}


To prove lemma \ref{lEstimatingThePhaseSpaceCommutatorInStrip}, we perform the same tricks of eliminating the $\qm$ localization to leave only $\pn$ localization. This will require the introduction of a smooth approximation to $\pni{x}$. 

\begin{proof}[of lemma \ref{lEstimatingThePhaseSpaceCommutatorInStrip}]
To prove this result, we need a smooth approximation to $\pni{x}$. Let $\pnibNOARG$ be a smooth, weakly increasing function with support on $[1/2,\infty)$ and identically one on $[1,\infty)$, and let $\pnicNOARG$ be a smooth, weakly decreasing function with support on $[0,2]$ and identically one on $[0,1]$. We extend these as even functions. Let 
\begin{align}
\pniaWITHlambda{\xi}=\pnibWITHlambda{\xi} \pnic{\xi} .
\hidenumber \end{align}

We first show $\pniaNOARG$ dominates $\pniWITHlambda{\bullet}$. Decomposing to harmonics and applying the spectral theorem in $\pn$, 
\begin{align}
\pniaWITHlambda{x} 
=& \pnibWITHlambda{x} \pnic{x} \hidenumber\\
\geq& \chi_{[\lambda^{-\delta},\infty)}(x) \chi_{[0,1]}(x) \hidenumber\\
\geq& \pniWITHlambda{x}, 
\label{ps16.2}
\end{align}
so that for $\psi\in L^2(\starman)$, 
\begin{align}
\| \pnia{\pn}\psi\|^2 \geq \|\pni{\pn}\psi\|^2 .
\hidenumber \end{align}
We also note
\begin{align}
\label{ps16.4}
\lambda^{-\delta} \pniaWITHlambda{\xi} \lesssim \xi\pnnc{\xi}\lesssim \xi\pnna{\xi}.
\end{align}

Since $\pnia{\pn}=\pnibNOARG(\pnplusdelta) \pnic{\pn}$, by the order reduction lemma, following the same argument which justifies equation \eqref{ps3.3}, if $F_1$ is a function with bounded derivative, such as $\qmnt{x}$, $\qmf{x}$, or $(-r_*^{-1}V_L')$, or a smooth compactly supported function, 

\begin{align}
[F_1(\qm),\pnia{\pn}]
=& [F_1(\qm),\pnibNOARG(\pnplusdelta)] \pnic{\pn} + \pnibNOARG(\pnplusdelta)[F_1(\qm),\pnic{\pn}] \hidenumber\\
=& L^{m+n+\delta-1}\bndd\pnic{\pn}+ \pnibNOARG(\pnplusdelta)L^{m+n-1}\bndd \hidenumber\\
=& L^{m+n+\delta-1}\bndd .
\hidenumber \end{align}

We start with the estimate involving the derivative operators. We relax the estimate by a factor of $\qmn{\qn}^\frac12$ and write $\dr$ in terms of $\pn$ before beginning to commute the localization, which generates the standard lower order terms. 
\begin{align}
\langle\psi,[H_1+H_3,\Gammann]\psi\rangle
\gtrsim& \|L^{n-\frac{\epsilon}{2}}\qmnt{\qn}^\frac12\dr\pnnan\chia\psi\|^2 +\ErrorTermsnPsi \hidenumber\\
\gtrsim& \|L^{1-\frac{\epsilon}{2}}\qmnt{\qn}\pn\pnnan\chia\psi\|^2 +\ErrorTermsnPsi \hidenumber\\
\gtrsim& \|L^{1-\frac{\epsilon}{2}}\pn\pnnan\qmnt{\qn}\chia\psi\|^2 +\ErrorTermsnPsi.
\hidenumber \end{align}
We now use equation \eqref{ps16.4} to introduce $\pnia{\pn}$ and continue commuting localizations. 
\begin{align}
\langle\psi,[H_1+H_3,\Gammann]\psi\rangle
\gtrsim& \|L^{1-\frac{\epsilon}{2}-\delta}\pnia{\pn}\qmnt{\qn}\chia\psi\|^2 +\ErrorTermsnPsi\hidenumber\\
\gtrsim& \|L^{1-\frac{\epsilon}{2}-\delta}\qmnt{\qn}\pnia{\pn}\chia\psi\|^2 +\ErrorTermsnPsi 
\label{ps17.3}
\end{align}
The last commutator is of order $2n+\delta-1$ in $L$. Since the introduction of $\pniaNOARG$ reduced the order of the entire expression to $1-\epsilon/2-\delta$, the term from the commutator is still lower order, of order $2n\leq1/2+n$. 

Turning to the potential term in equation \eqref{ps5.6}, we have 
\begin{align}
\langle\psi,[H_1+H_3,\Gammann]\psi\rangle
\gtrsim&\| L^{1-\frac{\epsilon}{2}} \qmf{\qn} (-r_*^{-1}V_L')^\frac12 \pnnan\chia\psi\|^2 + \ErrorTermsnPsi .
\hidenumber \end{align}
Using the same argument as in \eqref{enUsefToEliminateVLp}, we can eliminate the factor of $(-r_*^{-1}V_L')^\frac12$. 
\begin{align}
\langle\psi,[H_1+H_3,\Gammann]\psi\rangle
\gtrsim&\| L^{1-\frac{\epsilon}{2}} \qmf{\qn} \pnnan\chia\psi\|^2 + \ErrorTermsnPsi .
\hidenumber \end{align}

To continue, we drop a factor of $L^\delta$ and then start to commute localizations, and use $\pnna{x}\gtrsim\pniaWITHlambda{x}$. We need to drop the factor of $\delta$ powers of $L$ to control commutators involving $\pnia{\pn}$. 
\begin{align}
\langle\psi,[H_1+H_3,\Gammann]\psi\rangle
\gtrsim& \|L^{1-\frac{\epsilon}{2}-\delta} \qmf{\qn}\pnnan\chia\psi\|^2 + \ErrorTermsnPsi \hidenumber\\
\gtrsim& \|L^{1-\frac{\epsilon}{2}-\delta}\pnnan \qmf{\qn}\chia\psi\|^2 + \ErrorTermsnPsi \hidenumber\\
\gtrsim& \|L^{1-\frac{\epsilon}{2}-\delta}\pnia{\pn} \qmf{\qn}\chia\psi\|^2 + \ErrorTermsnPsi \hidenumber\\
\gtrsim& \|L^{1-\frac{\epsilon}{2}-\delta}\qmf{\qn} \pnia{\pn}\chia\psi\|^2 + \ErrorTermsnPsi .
\label{ps18.1}
\end{align}
We use \eqref{ps3.2} to combine \eqref{ps17.3} and \eqref{ps18.1}, and use \eqref{ps5.5} to control the contribution from $H_2$. 
\begin{align}
\langle\psi,[\linH,\Gammann]\psi\rangle
\gtrsim& \| L^{1-\frac{\epsilon}{2}-\delta} \pnia{\pn}\chia\psi\|^2 + \ErrorTermsnPsi .
\hidenumber \end{align}
We conclude by using \eqref{ps16.2} to replace $\pnia{\pn}$ by $\pni{\pn}$, 
\begin{align}
\langle\psi,[\linH,\Gammann]\psi\rangle
\gtrsim& 
\| L^{1-\frac{\epsilon}{2}-\delta} \pni{\pn} \chia\psi\|^2 + \ErrorTermsnPsi .
\hidenumber \end{align}
\end{proof}

\vspace{.5in}
{\bf\large Acknowledgement}

We would like to thank J. Sterbenz for providing valuable discussions concerning the results in \cite{BSterbenz}.

\end{document}